\documentclass[reqno, 11pt]{amsart}
\usepackage{amsfonts}
\usepackage[latin1]{inputenc} 
\usepackage[dvips]{color}

\newtheorem{theorem}{Theorem}[section]
\newtheorem{proposition}[theorem]{Proposition}
\newtheorem{lemma}[theorem]{Lemma}
\newtheorem{corollary}[theorem]{Corollary}

\newtheorem{remark}{Remark}[section]
\newtheorem{conjecture}[theorem]{Conjecture}

\newcommand{\nn}{\nonumber}
\newcommand{\oo}{\infty}

\newcommand{\ra}{\rightarrow}
\newcommand{\comp}{\circ}

\newcommand{\N}{\mathbb{N}}

\makeatletter
\@addtoreset{equation}{section}
\makeatother

\renewcommand\thefigure{\thesection.\@arabic\c@figure}
\renewcommand\thetable{\thesection.\@arabic\c@table}

\bibdata{prob}
\bibstyle{alpha}

\title[Spectrum for operators of stochastic machines]{Spectra of stochastic adding machines based on Cantor Systems of numeration}

\author{Ali Messaoudi$^1$, Glauco Valle$^2$}
\thanks{1. Supported by CNPq grant 305939/2009-2 and Universal CNPq project 482519/2012-6}
\thanks{2. Supported by CNPq grant 304593/2012-5 and Universal CNPq project 482519/2012-6}

\date{\today}

\address{
\newline
Ali Messaoudi
\newline
UNESP - Departamento de matem\'atica do Instituto de Bioci\^encias Letras e Ci\^encias Exatas.
\newline Rua Crist\'ov\~ao Colombo, 2265, Jardim Nazareth,
15054-000 - S\~ao Jos\'e do Rio Preto, SP, Brasil.
\newline
e-mail: {\rm \texttt{messaoud@ibilce.unesp.br}}
\newline
\newline
Glauco Valle
\newline
UFRJ - Departamento de m\'etodos estat\'{\i}sticos do Instituto de Matem\'atica.
\newline  Caixa Postal 68530, 21945-970, Rio de Janeiro, Brasil.
\newline
e-mail: {\rm \texttt{glauco.valle@im.ufrj.br}}
}


\begin{document}

\maketitle

\begin{abstract}
In this paper, we define a stochastic adding machine based on Cantor Systems of numeration. We also compute the parts of spectra of the transition operator associated to this stochastic adding machine in different Banach spaces as $c_0, c$ and $l_{\alpha},\; 1 \leq \alpha \leq +\infty$.
These spectra are connected to fibered Julia sets.

\end{abstract}

\section{Introduction}
\label{sec:intro}

Adding one to a non-negative integer $n$ can be performed through an algorithm that changes digits one by one in the expansion of $n$ on some system of numeration. Stochastic adding machines are time-homogeneous Markov Chains on the non-negative integers that are built based on a stochastic rule that prevents the algorithm to finish. Killeen and Taylor \cite{kt} have introduced this concept considering dyadic expansions. They found out and studied in detail an interesting relation between complex dynamics and the spectral properties of the transition operators of a particular class of dyadic stochastic machines: the spectra of the transition operators are the filled Julia sets of degree two polynomials. Thereafter stochastic machines have been studied on several other systems of numerations, see \cite{am,ms,msv}, generating a rich class of examples connecting probability theory, operator theory and complex dynamics. Before we give more details on previous works and on our motivations, let us introduce the stochastic machines we study in this paper.

Here we are going to consider stochastic machines on general Cantor systems of numeration. Denote the set of non negative integers by $\mathbb{Z}_+ = \{0,1,2,3,...\}$. Let us fix a sequence $\bar{d} = (d_l)_{l \ge 0}$ of positive integers such that $d_0 = 1$ and $d_j \ge 2$ for $j\ge 1$. Set
$$
\Gamma = \Gamma_{\bar{d}} := \Big\{ (a_j)_{j=1}^{+\oo} : a_j \in \{0,...,d_j-1\} ,\; j \ge 1 , \ \sum_{j=1}^{+\oo} a_j < \infty \Big\} \, ,
$$
and
$$
q_j = d_0 \, d_1 \, ... \, d_j \, , \ j \ge 0 \, .
$$
It is known (\cite{HW}, \cite{KT1}) that there is a one to one map from $\mathbb{Z}_+$ to $\Gamma$ that associates to each $n$ a sequence $(a_{j}(n))_{j=1}^{+\oo}$ such that
$$
n = \sum_{j=1}^{+\oo} a_{j}(n) q_{j-1} \, .
$$
The right hand side of the previous equality is called the q-expansion of $n$ and $a_{j}(n)$ is called the j th digit of the expansion. The map $n \mapsto n+1$ operates on $\Gamma$ in the following way: we define the counter $\zeta_n = \zeta_{\bar{d},n} := \min\{j\ge 1:a_{j}(n) \neq d_j-1\} $ then
$$
a_{j}(n+1) = \left\{
\begin{array}{cl}
0 &, \ j< \zeta_{n} \, , \\
a_{j}(n) + 1 &, \ j = \zeta_{n} \, , \\
a_{j}(n) &, \ j>\zeta_{n} \, .
\end{array}
\right.
$$
So an adding machine algorithm, that maps $n$ to $n+1$ using d-adic expansions by changing one digit on each step, is performed in $\zeta_{n}$ steps in the following way: the first $\zeta_{n}-1$ digits are replaced by zero recursively and in $\zeta_{n}th$ step we add one to the $\zeta_{n}th$ digit (basically we are adding one modulus $d_l$ on each step $l$). Note that $0\le \zeta_{n} \le l$ if $n \in [q_{l-1},q_l]$.



\medskip

Fix a sequence of strictly positive probabilities $\bar{p} = (p_j)_{j=1}^{+\oo}$. Suppose that at the j-th step of the adding machine algorithm, independently of any other step, the information about the counter get lost, thus making the algorithm to stop. This implies that the outcome of the adding machine is a random variable. We call this procedure the adding machine algorithm with fallible counter, or simply $\rm{AMFC}_{\bar{d},\bar{p}}$.

Formally, we fix a sequence $(\xi_j)_{j=1}^{+\oo}$ of independent random variables such that $\xi_j$ is a Bernoulli distribution with parameter $p_j \in (0,1]$. Define the random time $\tau = \inf\{ j : \xi_j = 0 \}$. Then the $\rm{AMFC}_{\bar{d},\bar{p}}$ is defined by applying the adding machine algorithm to $n$ and stopping at the step $\tau \wedge \zeta_{n}$ (this means that steps $j\ge \tau$ are not performed when $\tau < \zeta_{n}$).



\smallskip

Now fix a initial, possibly random, state $X(0) \in \mathbb{Z}_+$. We apply recursively the $\rm{AMFC}_{\bar{d},\bar{p}}$ to its successive outcomes starting at $X(0)$ and using independent sequences of Bernoulli random variables at different times. These random sequences are associated to the same fixed sequence of probabilities $(p_j)_{j=1}^{+\oo}$. In this way, we generate a discrete time-homogeneous Markov chain $(X(t))_{t \ge 0}$ which we call the $\rm{AMFC}_{\bar{d},\bar{p}}$ stochastic machine.

\smallskip

In \cite{kt}, the case $d_j = 2$, $p_j = p \in (0,1)$ is studied. Among other things, the authors show that the spectrum of the transition operator of the stochastic machine acting on $l^{\infty}$ is the filled-in Julia set of the degree two polynomial
$$
\frac{z^2 -(1-p)}{p} \, .
$$
Further spectral properties of the same transition operators and their dual acting on $c_0$, $c$, $l^\alpha$, $\alpha \ge 1$, are considered by Abdalaoui and Messaoudi in \cite{am}.

\smallskip

In \cite{msv}, Messaoudi, Sester and Valle have introduced the stochastic machines associated to non constant sequences $\bar{p}$, but $\bar{d} \equiv d$ constant. It is shown that the spectrum of its transition operator acting on $l^{\infty}$ is equal to the filled-in fibered Julia set $E$ defined by
$$
E_{d,\bar{p}} := \Big\{ z \in \mathbb{C} : \limsup_{j\ra +\oo} |\tilde{f}_{j}(z)| < +\oo \Big\}, \, .
$$
where
$\tilde{f}_{j} := f_{j} \comp ... \comp f_{1}$ for all $j \geq 1$ and
 $f_{j}:  \mathbb{C} \ra \mathbb{C}$, is the function defined by
$$
f_{j}(z) := \left( \frac{z-(1-p_j)}{p_j} \right)^d \, .
$$
It is also studied the topological properties of the filled-in fibered Julia sets $E_{d,\bar{p}}$. It is given sufficient conditions on the sequence $(p_n)_{n \geq 1}$ to ensure that $E_{d,\bar{p}}$ is a connected set, or has a finite or infinite number of connected components and also to ensure that the fibered Julia set $\partial E_{d,\bar{p}}$ is a quasi-circle.

\medskip

Motivated by \cite{am} and \cite{msv}, the aim of this paper is the analysis of the spectra of the transition operators of the more general $\rm{AMFC}_{\bar{d},\bar{p}}$ stochastic machines acting in $l^{\infty}$, $c_0$, $c$, $l^\alpha$, $\alpha \ge 1$. The decomposition of the spectra in their point, residual and continuous parts now depends on the particular choices of the sequences $\bar{d}$ and $\bar{p}$. But even in this far more general scenario we are able to give an almost complete description of the spectra answering positively a conjecture of Abdalaoui and Messaoudi that completely describes the residual spectrum in $l^1$.

In order to describe our results, we need to introduce some notation. Denote by $\sigma(\Omega,\bar{d},\bar{p})$, $\sigma_p(\Omega,\bar{d},\bar{p})$, $\sigma_r(\Omega,\bar{d},\bar{p})$ and $\sigma_c(\Omega,\bar{d},\bar{p})$ respectively the spectrum, point spectrum, residual spectrum and continuous spectrum of the transition operator of the $\rm{AMFC}_{\bar{d},\bar{p}}$ Markov chain acting as a linear operator on $\Omega \in \{ c_0, \ c, \ l^\alpha, \ 1 \le \alpha \le \oo \}$.  We obtain the following results:

\begin{itemize}
\item
The $\rm{AMFC}_{\bar{d},\bar{p}}$ Markov chain is null recurrent if and only if
$\prod_{j=1}^{+\oo} p_j = 0$,
otherwise the chain is transient.

\item

The spectrum of $S$ acting on $\Omega \in \{ c_0, \ c, \ l^\alpha, \ 1 \le \alpha \le \oo \}$ is equal to to the fibered filled  Julia set $ E_{\bar{d},\bar{p}} := \Big\{ z \in \mathbb{C} : \limsup_{j\ra +\oo} |\tilde{f}_{j}(z)| < +\oo \Big\}$ where $\tilde{f}_{j} := f_{j} \comp ... \comp f_{1}$ for all $j \geq 1$ and
 $f_{j}:  \mathbb{C} \ra \mathbb{C}$, is the function defined by
$
f_{j}(z) := \left( \frac{z-(1-p_j)}{p_j} \right)^{d_{j}}$.

\item
In $l^{\infty},$ the spectrum of $S$ is equal to the point spectrum $\sigma_p(l^\oo,\bar{d},\bar{p}) $ and hence the residual and continuous spectra are empty sets.
\item
In $\Omega \in \{ c_0, \ c, \ l^\alpha, \ 1 < \alpha \le \oo \}$, the residual spectrum is empty set.
\item
In $c_0$, the point spectrum is not empty if and only if $\lim p_n=1$, and in this case the $\sigma_p(c_0,\bar{d},\bar{p}) $ equals to the connected component of the interior of  $E_{\bar{d},\bar{p}} $ that contains $0$.
 Moreover, if $p_j \geq 2 (\sqrt{2}-1)$ for all integer $j \geq 1$ and $\bar{d}$ is constant, then  $\sigma_p(c_0,\bar{d},\bar{p})$ equals to the interior of  $E_{\bar{d},\bar{p}}$.
\item
In $c$, the point spectrum is equal to $\sigma_p(c_0,\bar{d},\bar{p})  \cup \{1\}$.
\item
In $l^{\alpha},\; \alpha  \geq 1$, if $\sum_{j=1}^{+\oo} (1-p_j)^\alpha$ diverges, then the point spectrum is empty. If $\bar{p}$ is monotone increasing and $\sum_{j=1}^{+\oo} (1-p_j)^\alpha$ converges, then $\sigma_p(l^{\alpha},\bar{d},\bar{p})= \sigma_p(c_0,\bar{d},\bar{p}) $.
 Note that, here we also consider the case $\alpha = 1$, where the convergence of $\sum_{j=1}^{+\oo} (1-p_j)$ is equivalent to $\prod_{j=1}^{+\oo} p_j > 0$.
\item
In $l^1$, if $(d_n)_{n \geq 0}$ is bounded and $\lim sup \; p_n <1$, then  the residual spectrum is contained in the boundary of $E_{\bar{d},\bar{p}} $ and equals to the countable set $X= \bigcup_{n=0}^{+\infty} f^{-n} \{1\} \setminus \bigcup_{n=0}^{+\infty} f^{-n} \{0\}$.

If  $\prod_{j=1}^{+\oo} p_j = 0$ (null-recurrent case), then  the residual spectrum  contains  the countable set $\bigcup_{n=0}^{+\infty} f^{-n} \{1\} \setminus \bigcup_{n=0}^{+\infty} f^{-n} \{0\}$.  In this case, we conjecture that $\sigma_r(l^1,\bar{d},\bar{p})= X$.
If $\prod_{j=1}^{+\oo} p_j > 0$ (transient case), $\sigma_r(l^1,\bar{d},\bar{p})\cap X= \emptyset$. In this case, we conjecture that $\sigma_r(l^1,\bar{d},\bar{p})= \emptyset $.

\end{itemize}

The paper is organized in the following manner: In the next section we shall define the $\rm{AMFC}_{\bar{d}}$ Markov chain and study when it is recurrent or transient.
In section three, we study the spectra of transition operator of the $\rm{AMFC}_{\bar{d}}$ chain  in $l^{\infty}$. In section four we focus the study in other Banach spaces as  $c_0, \ c $ and $ \ l^\alpha, \ 1 \le \alpha \le \oo$. The last section is devoted to the proofs of all technical lemmas.

\bigskip

\section{Transition operators and recurrence of $\rm{AMFC}_{\bar{d}}$ chains}
\label{sec:chain}

\smallskip

In this section $(X(t))_{t \ge 0}$ is an $\rm{AMFC}_{\bar{d},\bar{p}}$ stochastic machine associated to a sequence of non negative integers $\bar{d} = (d_j)_{j=1}^{+\oo}$, $d_j >1$, and of probabilities $\bar{p} = (p_j)_{j=1}^{+\oo}$, $p_j \in (0,1]$. This Markov chain is irreducible if and only if $p_j < 1$ for infinitely many j's. Moreover, when $p_j =1$ for every $j \ge 1$, we have that the $\rm{AMFC}_{\bar{d},\bar{p}}$ stochastic machine is the deterministic shift map $n \mapsto n+1$ on $\mathbb{Z}_+$.

Our first aim is to describe the transition probabilities of $(X(t))_{t \ge 0}$ which we denote $s(n,m) = s_{\bar{p},\bar{d}}(n,m):=P(X(t+1)=m|X(t)=n)$. They can be obtained directly from description of the chain: For every $n\ge 0$
\begin{equation}
\label{tp1}
s(n,m) =
\left\{
\begin{array}{cl}
(1-p_{r+1}) \prod_{j=1}^r p_j &, \ m = n - \sum_{j=1}^r (d_j-1) q_{j-1}  \, , \\
& \ \ r \le \zeta_n -1 \, , \ \zeta_n \ge 2 \, , \\
1-p_1 &, \ m=n \, , \\
\prod_{j=1}^{\zeta_n} p_j &, \ m=n+1 \, , \\
0 &, \ \textrm{otherwise} \, .
\end{array}
\right.
\end{equation}
From the exact expressions above, it is clear the self-similarity of the transition probabilities. Indeed, it is straightforward to verify that for every $j\ge 2$, we have
$$
s(n,m) =
\left\{
\begin{array}{cl}
s(n - a_{j-1}(n) q_{j-1}, m - a_{j-1}(n) q_{j-1}) &, \ q_{j-1} \le m \le q_{j}-1 \, , \\
0 &, \ \textrm{otherwise} \, ,
\end{array}
\right.
$$
for all $q_{j-1} \le n \le q_{j}-2$ (note that $\zeta_n \le j-1$ for this choice of $n$). Moreover, if $n=q_j-1$, we have $\zeta_n = j$, $s(q_j -1, q_j) = \prod_{l=1}^{j+1} p_l$ and
\begin{equation}
\label{tp4}
s(q_j -1, q_j - q_r) = (1-p_{r+1}) \prod_{l=1}^{r} p_l \, , \ 1\le r \le j \, .
\end{equation}
With the transition probabilities, we obtain the countable transition matrix of the $\rm{AMFC}_{\bar{d},\bar{p}}$ stochastic machine $S = S_{\bar{d},\bar{p}} = [s(n,m)]_{n,m \ge 0}$.

Note that $S$ is doubly stochastic if and only if $\prod_{j=1}^{+\oo} p_j = 0$. In fact $S$ is stochastic and the sum of coefficients of every column is $1$, except the first one whose sum is $1- \prod_{j=1}^{+\oo} p_j.$

\medskip

In the next proposition, we obtain a sufficient and necessary condition for recurrence of the $\rm{AMFC}_{\bar{d},\bar{p}}$ Markov chain.

\medskip

\begin{proposition}
\label{prop:recurrence}
The $\rm{AMFC}_{\bar{d},\bar{p}}$ Markov chain is null recurrent if and only if
\begin{equation}
\label{reccondition}
\prod_{j=1}^{+\oo} p_j = 0 \, .
\end{equation}
Otherwise the chain is transient.
\end{proposition}

\smallskip

The proof of Proposition \ref{prop:recurrence} is analogous to the proof for the case $\bar{d}$ constant that is found in \cite{msv}. We present it here for the sake of completeness.

\bigskip

\noindent \textbf{Proof:}  We start showing that condition (\ref{reccondition}) is necessary and sufficient to guarantee the recurrence of the $\rm{AMFC}_{\bar{d},\bar{p}}$ stochastic machine. From classical Markov chain Theory, the $\rm{AMFC}_{\bar{d},\bar{p}}$ stochastic machine is transient if and only if there exists a sequence $v=(v_j)_{j=1}^{+\oo}$ such that $0< v_j \le 1$ and
\begin{equation}
\label{rectran1}
v_j = \sum_{m = 1}^{+\infty} s(j,m) \, v_m \, , \ j \ge 1 \, ,
\end{equation}
i.e, $\tilde{S} \, v = v$ where $\tilde{S}$ is obtained from $S$ removing its first line and column. Indeed, in the transient case a solution is obtained by taking $v_m$ as the probability that $0$ is never visited by the $\rm{AMFC}_d$ Markov chain given that the chain starts at state $m$ (see the discussion on pages 42-43 of Chapter 2 in \cite{l} and also \cite{O}).

\smallskip

Suppose that $v=(v_j)_{j=1}^{+\oo}$ satisfies the above conditions. We claim that
\begin{equation}
\label{rectran2}
v_{q_l + j} = v_{q_l} \, , \, \textrm{ for every } \, l\ge 0 \, \textrm{ and } \, j \in \{ 1, ..., (d_{l+1}-1)q_l-1 \} \, .
\end{equation}
The proof follows from induction. Indeed, for $j \in \{ 1, ..., (d-1)q_l-1 \}$, suppose that $v_{q_l} = v_{q_l + r}$, for all $0 \le r \le j-1$  we have that
$v_{q_l+j-1}  =  \sum_{m=1}^j s(q_l+j-1,m) \, v_{m}$. Since
$s(q_l+j-1,m)=0$  for all $0 \leq m < q_l$, we have
\begin{eqnarray}
\label{rectran3}
v_{q_l+j-1} & = & \sum_{r=0}^j s(q_l+j-1,q_l+r) \, v_{q_l + r} \nn \\
& = & \left( \sum_{r=0}^j s(q_l+j-1,q_l+r) \right) \, v_{q_l} \nn \\
& & \quad \quad + s(q_l+j-1,q_l+j) (v_{q_l+j} - v_{q_l})  \, .
\end{eqnarray}
Using the fact that $j \in \{ 1, ..., (d-1)q_l-1 \}$ note that
$$
\sum_{r=0}^j s(q_l+j-1,q_l+r) = 1.
$$
Thus, since $s(q_l+j-1,q_l+j)>0$, from (\ref{rectran3}), we have that $v_{q_l+j} = v_{q_l+j-1} = v_{q_l}$. This proves the claim.

\smallskip

It remains to obtain $v_{q_{l+1}}$ from $v_{q_{l}}$ for $l\ge 0$. First note that (\ref{rectran2}) implies $v_{q_{l+1} - q_r} = v_{q_l}$ for $0\le r \le l$. From the transition probabilities expression in (\ref{tp4}), if we put $p_0=1$, we have that
\begin{eqnarray*}
v_{q_l} & = & v_{q_{l+1}-1} \ = \ (p_0...p_{l+2}) v_{q_{l+1}} + \sum_{r=0}^{l} (p_0...p_r - p_0...p_{r+1}) v_{q_{l+1} - q_r} \, . \\
& = & (p_0...p_{l+2}) v_{q_{l+1}} + (1 - p_0...p_{l+1}) v_{q_{l}}
\end{eqnarray*}
Therefore for every $l\ge 1$
$$
v_{q_{l}} = \frac{v_{q_{l-1}}}{p_{l+1}} = \frac{v_1}{\prod_{j=2}^{l+1} p_j} \, .
$$
From this equality, we conclude  that $v$ exists and the chain is transient if and only if
$$
\prod_{j=1}^{+\oo} p_j > 0 \, .
$$

\smallskip

Now suppose that we are in the recurrent case. Since $S$ is a irreducible countable doubly stochastic matrix, it is simple to verify that the $\rm{AMFC}_{\bar{d},\bar{p}}$ have no finite invariant measure and then cannot be positive recurrent.
$\square$

\bigskip

\section{$l^\infty$ Spectra of transition operators of $\rm{AMFC}_{\bar{d},\bar{p}}$ chains}
\label{sec:spectra}

\smallskip

In this section we discuss the spectra of transition operators of $\rm{AMFC}_{\bar{d},\bar{p}}$ stochastic machines. We consider the usual notation for the Banach spaces $(l^\alpha,\| \cdot \|_\alpha)$, $1 \le \alpha \le \oo$, $(c_0,\| \cdot \|_\oo)$ and $(c,\| \cdot \|_\oo)$, i.e.: for $w = (w(n))_{n \ge 0} \in \mathbb{C}^{\mathbb{Z}_+}$, we have
$$
\| w \|_\oo =  \sup_{n \ge 0} |w(n)| < \oo \, , \quad
\| w \|_\alpha = \Big( \sum_{n\ge 0} |w(n)|^\alpha \Big)^\frac{1}{\alpha} \, , \ 1 \le \alpha < \oo \, ,
$$
and
$$
l^\oo = l^\oo(\mathbb{Z}_+) = \{ w \in \mathbb{C}^{\mathbb{Z}_+} : \| w \|_\oo < \oo \} \, ,
$$
$$
l^\alpha = l^\alpha(\mathbb{Z}_+) = \{ w \in \mathbb{C}^{\mathbb{Z}_+} : \| w \|_\alpha < \oo \} \, ,
$$
$$
c = c(\mathbb{Z}_+) = \{ w \in l^\oo : w \textrm{ is convergent}  \} \, ,
$$
$$
c_0 = c_0 (\mathbb{Z}_+) = \{ w \in c : \lim_{n\rightarrow \oo} w(n) = 0 \} \, .
$$
The transition operator of the $\rm{AMFC}_{\bar{d},\bar{p}}$ stochastic machine is the bounded linear operator on $l^\oo$ induced by $S_{\bar{d},\bar{p}}$ as
$$
[S_{\bar{d},\bar{p}} w](n) = \sum_{m=0}^\infty s(n,m) w(m) \, , \ n \in \mathbb{Z}_+ \, ,
$$
for every $w \in l^\oo$. We also denote the transition operator by $S_{\bar{d},\bar{p}}$. The dual transition operator of the $\rm{AMFC}_{\bar{d},\bar{p}}$ stochastic machine, which we denote by $S^\prime_{\bar{d},\bar{p}} \, ,$ is the bounded linear operator on $l^\oo$ induced by $S_{\bar{d},\bar{p}}$ as
$$
[S^\prime_{\bar{d},\bar{p}} w](n) = [w S_{\bar{d},\bar{p}}](n) = \sum_{m=0}^\infty w(m) s(m,n) \, , \ n \in \mathbb{Z}_+ \, ,
$$
for every $w \in l^\oo$.

The matrix $S_{\bar{d},\bar{p}}$ is stochastic with columns also having its sums uniformly bounded above by one. Therefore we can show, analogously to Proposition 4.1 in \cite{am}, that the restrictions of $S_{\bar{d},\bar{p}}$ and $S^\prime_{\bar{d},\bar{p}}$ to $\Omega \in \{c_0, c, l^\alpha, 1\le \alpha < \oo\}$ are well defined bounded linear operators on $\Omega$.

We call $S^\prime_{\bar{d},\bar{p}}$ the dual transition operator in order to simplify the text. From usual operator theory, the dual of $S_{\bar{d},\bar{p}}$ acting on $c_0$, $c$ is $S^\prime_{\bar{d},\bar{p}}$ acting on $l^1$ and the dual of $S_{\bar{d},\bar{p}}$ acting on $l^\alpha,\; 1 \leq \alpha < \infty$ is $S^\prime_{\bar{d},\bar{p}}$ acting on $l^\frac{\alpha}{\alpha -1}$.

\medskip

We recall that given a complex Banach space $E$  and $T: E
\rightarrow E$ a continuous linear operator, the spectrum of the
operator $T$ can be partitioned into three subsets (see for instance \cite{Y}):

\begin{enumerate}

\item The point spectrum $ \sigma_{p}(T) = \{\lambda \in \mathbb{C} :
T - \lambda I  \mbox{ is not injective} \}$.

\item The continuous spectrum
$\sigma_c(T) = \{\lambda \in \mathbb{C}: T - \lambda I \
\mbox{ is injective }, \\
 \overline{(T - \lambda I)E} = E, \ (T -
\lambda I)E \neq E \}$, where $\overline{(T - \lambda I)E}$ is the closure of $(T - \lambda I)E$ in $E$.

\item The residual spectrum
$\sigma_r(T) = \sigma(T) \setminus (\sigma_{p}(T) \cup \sigma_c(T)) = \{\lambda \in \mathbb{C}:
\lambda I - T \ \mbox{is injective}, \ \overline{(\lambda I - T)E}
\neq E \}$.
\end{enumerate}

\medskip

In order to describe the spectrum of $S_{\bar{d},\bar{p}}$, we need to introduce more notation.
We denote by  $\mathbb{D}(w,r) = \{ z \in \mathbb{C} : |w-z| < r \}$ and $\overline{\mathbb{D}(w,r)} = \{ z \in \mathbb{C} : |w-z| \leq r \}$. Let $f_{j}:  \mathbb{C} \ra \mathbb{C}$, $j \ge 1$, be the function defined by
$$
f_{j}(z) := \left( \frac{z-(1-p_j)}{p_j} \right)^{d_j} \, .
$$
Also set $\tilde{f}_0$ as the identity function on $\mathbb{C}$, $\tilde{f}_{j} := f_{j} \comp ... \comp f_{1}$, $j \ge 1$, and
$$
E_{\bar{d},\bar{p}} := \Big\{ z \in \mathbb{C} : \limsup_{j\ra +\oo} |\tilde{f}_{j}(z)| < +\oo \Big\} \, .
$$

\medskip

\begin{lemma}\label{bound}
The set $E_{\bar{d},\bar{p}}$ is included in $\overline{\mathbb{D}(1-p_1,p_1)}$. Moreover, for all $z\in E_{\bar{d},\bar{p}}$ and $j\geq 1$, $\tilde{f}_{j}(z)$ belongs to the disk $\overline{\mathbb{D}(1-p_{j+1},p_{j+1})}$.
\end{lemma}

\medskip

\begin{corollary}
\label{cor:bound}
We have that
\begin{eqnarray*}
E_{\bar{d},\bar{p}} & = & \overline{\mathbb{D}(1-p_1,p_1)} \cap \bigcap_{j=1}^\infty \tilde{f}_j^{-1} \big( \overline{\mathbb{D}(1-p_{j+1},p_{j+1})} \big) \\
& = &
\overline{\mathbb{D}(0,1)} \cap \bigcap_{j=1}^\infty \tilde{f}_j^{-1} \big( \overline{\mathbb{D}(0,1) \big)} \, .
\end{eqnarray*}
\end{corollary}

\bigskip

Let $\Omega \in \{ c_0, \ c, \ l^\alpha, \ 1 \le \alpha \le \oo \}$.
Denote by $\sigma(\Omega,\bar{d},\bar{p})$, $\sigma_p(\Omega,\bar{d},\bar{p})$, $\sigma_r(\Omega,\bar{d},\bar{p})$ and $\sigma_c(\Omega,\bar{d},\bar{p})$ respectively the spectrum, point spectrum, residual spectrum and continuous spectrum of $S_{\bar{d},\bar{p}}$ acting as a linear operator on $\Omega$. We replace $\sigma$ by $\sigma^\prime$ to represent the spectrum and its decomposition for the dual transition operator $S^\prime_{\bar{d},\bar{p}}$.

\medskip

Let us start with the analysis of a simple and well-known but important case. Recall that for $p_j = 1$, $j\ge 1$, we have that the $\rm{AMFC}_{\bar{d},\bar{p}}$ stochastic machine is the deterministic map $n \mapsto n+1$ on $\mathbb{Z}_+$ for any sequence $\bar{d}$. In this case, $\sigma(l^\oo,\bar{d},\bar{p}) = \sigma_p(l^\oo,\bar{d},\bar{p}) = \overline{\mathbb{D}(0,1)}$. Since, $\tilde{f}_{j} (z) = z^j$, $z \in \mathbb{C}$, we have that $(\tilde{f}_{j} (z))_{j \ge 1}$ is bounded, if and only if, $z \in \overline{\mathbb{D}(0,1)}$. Therefore, $E_{\bar{d},\bar{p}} = \sigma(l^\oo,\bar{d},\bar{p}) = \sigma_p(l^\oo,\bar{d},\bar{p}) = \overline{\mathbb{D}(0,1)}$.

For the general case, we have:

\bigskip

\begin{theorem}
\label{teo:ps}
For every sequences $\bar{p} \in (0,1]^\mathbb{N}$ and $\bar{d}= (d_{i})_{i \geq 0}$ where $d_0=1$ and  $(d_{i})_{i \geq 1} \subset \{2,3,4,...\}^\mathbb{N}$,
$$
E_{\bar{d},\bar{p}} = \sigma(l^\oo,\bar{d},\bar{p}) = \sigma_p(l^\oo,\bar{d},\bar{p}) \, .
$$
\end{theorem}

\bigskip

The proof of Theorem \ref{teo:ps} follows directly from the next two propositions:

\bigskip

\begin{proposition}
\label{prop:ps}
The point spectrum of $S_{\bar{d},\bar{p}}$ in $l^\oo$ is equal to $E_{\bar{d},\bar{p}}$.
\end{proposition}

\medskip

\begin{proposition}
\label{spectruminjulia}
For any $\Omega \in \{ c_0, \ c, \ l^\alpha, \ 1 \le \alpha \le \oo \}$,
$\sigma(\Omega,\bar{d},\bar{p}) \subset E_{\bar{d},\bar{p}}$.
\end{proposition}

\bigskip

The rest of this section will be devoted to the proof of the previous propositions. To prove the proposition \ref{prop:ps}, we give in the next lemma an explicit characterization of the eigenvectors of $S_{\bar{d},\bar{p}}$.

\medskip

\begin{lemma}
\label{lemma:ps}
A sequence $v \in l^\oo$ is an eigenvector of $S_{\bar{d},\bar{p}}$ associated to an eigenvalue $\lambda$, if and only if, for some $v(0) \in \mathbb{C}$, $v = v(0) \cdot v_\lambda$ with $v_\lambda$ given by
\begin{equation}
\label{ps1}
v_\lambda (n) = \prod_{r=1}^{\infty} (\iota_{\lambda}(r))^{a_{r}(n)} \, , \ n \ge 0 \, ,
\end{equation}
where $a_{r}(n)$ is the $r$ digit of $n$ in its $q$-expansion and
\begin{equation}
\label{ps2}
\iota_{\lambda}(r) = (\, h_r \comp \tilde{f}_{r-1} \, ) (\lambda)
\end{equation}
for
\begin{equation}
\label{ps3}
h_r(z) = \frac{z}{p_r} - \frac{1-p_r}{p_r} \, .
\end{equation}
\end{lemma}

\bigskip

Note that in (\ref{ps1}) we only have a finite number of terms in the product that are distinct from 1.

The proofs of Proposition \ref{prop:ps} and Lemma \ref{lemma:ps} are analogous to the case $\bar{d}$ constant presented in \cite{msv}. We present both of them here for the sake of completeness.

\bigskip

\noindent \textbf{Proof of Proposition \ref{prop:ps}:}
Since $S_{\bar{d},\bar{p}}$ is stochastic, its spectrum is a subset of $\overline{D(0,1)}$. By Lemma \ref{lemma:ps}, $\lambda \in  \overline{D(0,1)}$ is eigenvalue of $S_{\bar{d},\bar{p}}$, if and only if, $v_\lambda \in l^\oo$. In this case, $v_\lambda$ is, up to multiplication by a constant, the unique eigenvector of $S_{\bar{d},\bar{p}}$ in $l^\oo$ associated to $\lambda$.

We are going to show that $(|\iota_{\lambda}(j)|)_{j=1}^{+\oo}$ is bounded above by one if and only if $\lambda \in E_{\bar{d},\bar{p}}$, otherwise it is unbounded. Therefore, from (\ref{ps1}), we have that $v_\lambda$ is a well defined element of $l^\oo$ if and only if $\lambda \in E_{\bar{d},\bar{p}}$. Thus Proposition \ref{prop:ps} holds.

If $\lambda \in E_{\bar{d},\bar{p}}$ then $\tilde{f}_{r-1}(\lambda)$ is uniformly bounded and according to Lemma \ref{bound}, for all $r\geq 1$
$$
\tilde{f}_{r-1}(\lambda)\in \overline{\mathbb{D} (1-p_{r},p_{r})}.
$$
Since $h_r$ maps  $\overline{\mathbb{D} (1-p_{r},p_{r})}$ on $\overline{\mathbb{D} (0,1)} $ we deduce that
$|\iota_{\lambda}(r)|\leq1 $ and $(|\iota_{\lambda}(r)|)_{r=1}^{+\oo}$ is bounded above by one.
Indeed, assume that there exists $j_0 \in \mathbb{N}$ such that $|\iota_{\lambda}(j_0)|>1$. Then $\lambda \not \in E_{\bar{d},\bar{p}}$. Thus  $|\tilde{f}_{r}(z)|>1$ for some $r>1$. We deduce, by (\ref{lowerbound2}) in the proof of Lemma \ref{bound}, that $\lim_{j \to +\infty} \vert  \tilde{f}_{j}(z) \vert= +\infty$. Hence $\lim_{j \to +\infty}|\iota_{\lambda}(j)|)= +\infty$.

Conversely, suppose $|\iota_{\lambda}(r)|\leq 1$ for all $r$. From (\ref{lowerbound1}) also in the proof of Lemma \ref{bound}, we know that for any $|z|>1$
$$
|h_r(z)| \ge |z|. \,
$$
Thus, if $|\tilde{f}_{r-1}(\lambda)| > 1$ for some $r > 0$, we have
$$
|\iota_{\lambda}(r)| = |h_r(\tilde{f}_{r-1}(\lambda))| \ge |\tilde{f}_{r-1}(\lambda)| > 1 \, ,
$$
which yields a contradiction to the fact that $|\iota_{\lambda}(r)|\leq 1$. Hence $|\tilde{f}_{r-1}(\lambda)|\leq 1$ for all $r$ and, by definition, $\lambda \in E$.

Hence $(|\iota_{\lambda}(j)|)_{j=1}^{+\oo}$ is bounded above by one if and only if $\lambda \in E_{\bar{d},\bar{p}}$.
$\square$

\bigskip \medskip

\noindent \textbf{Proof of Proposition \ref{spectruminjulia}:}
We prove here that $\sigma(\Omega,\bar{d},\bar{p}) \subset E_{\bar{d},\bar{p}}$. Here we denote by $\tau:\mathbb{Z}_+ \ra \mathbb{Z}_+$ the shift map $\tau(n)=n+1$ and by $\bar{p}_n := (p_{n+j})_{j=0}^\oo$, $\bar{d}_n := (d_{n+j})_{j=0}^\oo$.

\medskip

Put
$$
\tilde{S}_{\bar{p}} := \frac{S_{\bar{p}} - (1-p_1)I}{p_1} \, ,
$$
which is also a stochastic operator acting on $\mathbb{Z}_+$.  It is associated to a irreducible  Markov chain with period $d$. Thus $\tilde{S}_{\bar{d},\bar{p}}^{d_1}$ has $d_1$ communication classes.

\medskip

It is straightforward to verify that the communication classes of $\tilde{S}_{\bar{d},\bar{p}}^{d_1}$ are
$$
\{ \ j\in \mathbb{N} \ : \ j=n \mod d_1 \} \, , \quad 0\le n \le d_1-1 \, .
$$
Furthermore, $\tilde{S}_{\bar{d},\bar{p}}^{d_1}$ acts on each of these classes as a copy of $S_{\bar{d}_2,\bar{p}_2,}$. Therefore, the spectrum of $\tilde{S}_{\bar{d},\bar{p}}^{d_1}$ is equal to the spectrum of $S_{\bar{d}_2,\bar{p}_2}$. Since, $\tilde{S}_{\bar{d},\bar{p}}^{d_1} = \tilde{f}_{1}\big( S_{\bar{d},\bar{p}} \big)$, by the Spectral Mapping Theorem, we have that
$$
\tilde{f}_{1}\big( \sigma(\Omega,\bar{d},\bar{p}) \big) = \sigma(\Omega,\bar{d}_2,\bar{p}_2) \, .
$$
By induction, we have that
$$
\tilde{f}_{j+1}\big( \sigma(\Omega,\bar{d},\bar{p}) \big) = \sigma(\Omega,\bar{d}_{j+1},\bar{p}_{j+1}) \, ,
$$
for every $j\ge 1$. Since, $S_{\bar{d}_{j+1},\bar{p}_{j+1}}$ is a stochastic operator, its spectrum is a subset of $\overline{\mathbb{D}(0,1)}$. Therefore
$$
| \tilde{f}_{j+1}\big( \lambda \big) | \le 1 \, ,
$$
for every $j$ and $\lambda \in \sigma(\Omega,\bar{d},\bar{p})$. This implies that $\sigma(\Omega,\bar{d},\bar{p}) \subset E_{\bar{d},\bar{p}}$.
$\square$

\bigskip

\section{Spectra of transition operators of $\rm{AMFC}_{\bar{d},\bar{p}}$ chains on other Banach spaces}
\label{sec:spectra}

\medskip

By Proposition \ref{spectruminjulia}, for any $\Omega \in \{ c_0, \ c, \ l^\alpha, \ 1 \le \alpha \le \oo \}$, $\sigma(\Omega,\bar{d},\bar{p}) \subset E_{\bar{d},\bar{p}}$. We will indeed show that $\sigma(\Omega,\bar{d},\bar{p}) = E_{\bar{d},\bar{p}}$ as in the case $\Omega = l^\infty$. From this point we can ask ourselves about the decomposition of  $\sigma(\Omega,\bar{d},\bar{p})$ in its point, residual and continuous parts. We will see that this decomposition depends on the parameters of the $\rm{AMFC}_{\bar{d},\bar{p}}$ Markov Chains generating a rich class of examples.

\medskip

\begin{theorem}
\label{spectrumequaljulia}
For $\Omega \in \{ c_0, \ c, \ l^\alpha, \ 1 \le \alpha < \oo \}$, we have that $\sigma(\Omega,\bar{d},\bar{p}) = E_{\bar{d},\bar{p}}$.
\end{theorem}

\bigskip

The main step to prove Theorem \ref{spectrumequaljulia} is the following result:

\bigskip

\begin{lemma}
\label{lemma:spectrumequaljulia}
For $1 \le \alpha < \oo$, every $\lambda \in E_{\bar{d},\bar{p}} \setminus \sigma_p (l^\alpha,\bar{d},\bar{p})$ belongs to the approximate point spectrum of $S_{\bar{d},\bar{p}}$ acting on $l^\alpha$.
\end{lemma}

\bigskip

\noindent \textbf{Proof of Theorem \ref{spectrumequaljulia}:} Assume that $\Omega \in \{c_0,\ c \}$. Then, by duality and Phillips Theorem, we obtain
$$
\sigma(\Omega,\bar{d},\bar{p}) = \sigma^\prime (l^1,\bar{d},\bar{p})
= \sigma(l^\oo,\bar{d},\bar{p}) = E_{\bar{d},\bar{p}} \, .
$$
Now, assume $\Omega= l^{\alpha},\; 1 \le \alpha < \oo$. According to Proposition \ref{spectruminjulia}, it is enough to prove that $E_{\bar{d},\bar{p}} \subset \sigma(l^\alpha,\bar{d},\bar{p})$. This follows from Lemma \ref{lemma:spectrumequaljulia}, since every point in $E_{\bar{d},\bar{p}}$ is in the point or approximate point spectrum of $S_{\bar{d},\bar{p}}$.
$\square$

\bigskip

From now  we study the decomposition of $\sigma(\Omega,\bar{d},\bar{p})$, $\Omega \in \{ c_0, \ c, \ l^\alpha, \ 1 \le \alpha < \oo \}$ in its point, residual and continuous parts. We start with the case $\Omega \neq l^1$ and we consider the $l^1$ case later due to its particularities.

\medskip

\begin{proposition}
\label{prop:residualempty}
For  $\Omega \in \{ c_0, \ c, \ l^\alpha, \ \alpha > 1 \}$, we have that $\sigma_r(\Omega,\bar{d},\bar{p})$ is empty.
\end{proposition}

\medskip

The proof of Proposition \ref{prop:residualempty} relies on duality. In this direction a proper representation for the left eigenvectors of $S^\prime_{\bar{d},\bar{p}}$ is useful.

\medskip

\begin{lemma}
\label{lemma:left-eigenvector}
A sequence $v^\prime \in l^\oo$ is an eigenvector of $S^\prime_{\bar{d},\bar{p}}$ associated to an eigenvalue $\lambda$, if and only if, for some $v^\prime(0) \in \mathbb{C}$, $v^\prime = v(0) \cdot v^\prime_\lambda$ with $v^\prime_\lambda$ given by
\begin{equation}
\label{eq:rs1}
v^\prime_\lambda(m) = \Big( \prod_{r=1}^{\infty} (\iota_{\lambda}(r))^{a_{r}(m)} \Big)^{-1} = \frac{1}{v_\lambda (m)}  \, , \, \textrm{ for every } \, m \ge 1.
\end{equation}
where $\iota_{\lambda}(r)$ is defined as in statement of Lemma \ref{lemma:ps}.
\end{lemma}

\medskip

\begin{remark}
Since $\lambda v'_\lambda (0) = [ S^\prime_{\bar{d},\bar{p}} v^\prime_\lambda ](0)$ (see { \ref{eq:rs2}}), we have that
$$
(1- p_1 - \lambda)v'_\lambda (0) + \sum_{i=1}^{+\infty} \left((1-p_{i+1}) \prod_{j=1}^{i} p_j\right) v'_\lambda (q_i -1)=0
$$
which is equivalent to
\begin{eqnarray}
\label{tb}
\iota_{\lambda}(1) =  \sum_{i=1}^{+\infty}
\frac{(1-p_{i+1}) \prod_{j=2}^{i} p_j}{\prod_{r=1}^{i-1} (\iota_{\lambda}(r))^{d_r-1}}.
\end{eqnarray}
\end{remark}

\bigskip

\noindent \textbf{Proof of Proposition \ref{prop:residualempty}:}
Fix the space $\Omega \in \{ c_0, c, l^\alpha, \alpha > 1 \}$.  From classical operator theory, we have that the residual spectrum is a subset of point spectrum of the dual operator. Then, we have to prove that does not exist $\lambda \in E_{\bar{d},\bar{p}}$ and $w = (w_n)_{n\ge 1} \in l^1$ such that $w \neq 0$ and $S^\prime_{\bar{d},\bar{p}} w = \lambda w$.

For $\lambda \in E_{\bar{d},\bar{p}}$, $v_{\lambda}$ is uniformly bounded above by 1, then $v_{\lambda}$ is uniformly bounded below by 1. By Lemma \ref{lemma:left-eigenvector}, we see that if $S^\prime_{\bar{d},\bar{p}} w = \lambda w$, with $\lambda \in E_{\bar{d},\bar{p}}$, then $\vert w(m) \vert \ge \vert w(0) \vert \, \vert v^\prime_\lambda(m) \vert \ge \vert w(0) \vert $, for every $m$. Hence $v \in l^1$ only if $v \equiv 0$.
$\square$

\bigskip

\begin{proposition}
\label{ptspec}
For $\Omega \in  \{c_0,  l^\alpha, \alpha \geq 1\}$, if $\bar{p}$ does not converge to $1$, then $\sigma_{p} (\Omega,\bar{d},\bar{p})$ is  empty.
\end{proposition}

\medskip

\noindent \textbf{Proof:}
Fix the space $\Omega \in \{ c_0, \; l^\alpha, \; \alpha \geq 1 \}$.  Since $\Omega \subset l^\infty$, we have that $\lambda$ is an eigenvalue of $S_{\bar{d},\bar{p}}$ on $\Omega$ only if it is an eigenvector of $S_{\bar{d},\bar{p}}$ on $l^\infty$ having an eigenvector in $\Omega$. By Proposition \ref{prop:ps}, for this to happen is necessary that $\lim_{j \ra \infty} \iota_\lambda(j) = 0$ for some $\lambda \in E_{\bar{d},\bar{p}}$.
Since
\begin{eqnarray}
\label{xc}
\iota_{\lambda}(j+1)= \frac{\iota_{\lambda}(j)^ { d_{j} } }{p_{j+1}}- \frac{1- p_{j+1}}{p_{j+1}},\; \forall j \in \mathbb{Z}^{+} \, ,
\end{eqnarray}
if $\bar{p}$ does not converge to $1$, then $\iota_\lambda(j)$ does not converge to $0$. Thus  $\sigma_{p} (\Omega,\bar{d},\bar{p}) = \emptyset$.
$\square$

\bigskip

Let us point out that the condition on $\bar{p}$ in the statement is necessary. Recall that for the shift, $p_j = 1$ for all $j$, we have $\sigma_{p} (\Omega,\bar{d},\bar{p}) = \mathbb{D}(0,1)$ for $\Omega \in  \{c_0, , l^\alpha, \alpha \ge 1\}$. Indeed the condition on the statement is necessary even when $p_j \neq 1$ for every $j \ge 1$. This is shown in the next result.

\medskip

From Propositions \ref{prop:residualempty} and \ref{ptspec}, we have:

\medskip

\begin{theorem}
For $\Omega \in  \{c_0, l^\alpha, \alpha \geq 1\}$, if $\bar{p}$ does not converge to $1$, then $\sigma (\Omega,\bar{d},\bar{p}) = \sigma_{c} (\Omega,\bar{d},\bar{p})$.
\end{theorem}

\medskip

\subsection{Case $\Omega= c_0$ or $\Omega= c$}.

\smallskip

For all integer  $n \in \mathbb{N}$, let $g_n: \mathbb{C} \to \mathbb{C}$ be
the function defined by $g_n (\lambda)= i_{\lambda}(n)$ for all $\lambda \in \mathbb{C}$.

\medskip

\begin{proposition}
\label{lemma:ptnotempty}
If $\lim_{j \rightarrow \oo} p_j = 1$, then $int (\sigma_{p} (c_0,\bar{d},\bar{p}))$ is not empty.
Precisely, there exists  a real number $r >0$ and an integer $j_0 \geq 1$  such that for all integer $j \geq j_0$ the open set  $g_{j}^{-1} \big( B(0, r) \big) \subset \sigma_{p} (c_0,\bar{d},\bar{p})$.
\end{proposition}
\medskip


\noindent \textbf{Proof:}
Put $\rho = 2 (\sqrt{2} - 1)$. We are going to show the following assertion: If $\lim_{j \ra \oo} p_j = 1$ and there exists $j_0$ such that $\inf_{j\ge j_0} p_j \ge \rho$ and $|\iota_\lambda (j_0)| \le r := \rho / 2$, then $\lim_{j \ra \oo} |\iota_\lambda (j)| = 0$. Since $\lim_{j \ra \oo} |\iota_\lambda (j)| = 0$ implies that $\lambda \in  \sigma_{p} (c_0,\bar{d},\bar{p})$, we have that $g_{j}^{-1} \big( B(0, r) \big) \subset \sigma_{p} (c_0,\bar{d},\bar{p})$ for $j \ge j_0$.

To prove the previous assertion, we construct, for any fixed $\eta \in (1,2)$  a subsequence $(j_k)_{k\ge 0}$ such that
$$
 |\iota_\lambda(j)| \le r^{\eta^k} , \ \textrm{ for every } j_k \le j \le j_{k+1} -1 \, .
$$
Since $r<1$ and $\eta > 1$, the assertion holds.

We construct $(j_k)_{k\ge 0}$ by induction. For $k=0$, take $j \ge j_0$ and suppose $|\iota_\lambda(j)| \le r$, then
$$
|\iota_{\lambda}(j + 1)| \le \frac{|\iota_{\lambda}(j)|^{d_{j}}}{p_j} + \frac{1-p_j}{p_j} \le
\frac{r^2}{\rho} + \frac{1-\rho}{\rho} = r \, .
$$
Therefore, $|\iota_\lambda(j)| \le r$ for every $j \ge j_0$. Now fix $k \ge 1$ and suppose that there exists $j_k > j_0$ such that $|\iota_\lambda(j)| \le r^{\eta^k}$ for every $j \ge j_k$. Since $r^{2 \eta^k} \le r^{\eta^{k+1}}$ and $\lim p_j = 1$, there exists $j_{k+1} > j_k$ such that
$$
\frac{r^{2 \eta^k}}{\rho_{k+1}} + \frac{1-\rho_{k+1}}{\rho_{k+1}} \le r^{\eta^{k+1}} \, ,
$$
where $\rho_k = p_{j_k} = \inf_{j \ge j_{k+1}} p_j$. Then, for every $j \ge j_{k+1}$,
\begin{equation}
\label{eq:iconv}
|\iota_\lambda(j)| \le \frac{|\iota_\lambda(j)|^{d_{j}}}{p_j} + \frac{1-p_j}{p_j} \le
\frac{r^{2 \eta^k}}{\rho_{k+1}} + \frac{1-\rho_{k+1}}{\rho_{k+1}} \le r^{\eta^{k+1}} \, .
\end{equation}
Therefore $\lim_{j \ra \oo} |\iota_\lambda (j)| = 0$ and the proof is complete. $\square$

\vspace{1em}

\begin{proposition}
\label{composant}
If $\lim_{j \rightarrow \oo} p_j = 1$, then
 $\sigma_{p} (c_0,\bar{d},\bar{p})$ equals to the connected component  of  $ int (E_{\bar{d},\bar{p}}) $ that contains $0$.
\end{proposition}

\noindent \textbf{Proof:}
Let $V$ be the connected component of $int (E_{\bar{d},\bar{p}}) $ that contains $0$.
Let $O = B(0,r)$ be a neighborhood of $0$ where $r$ is as in Proposition \ref{lemma:ptnotempty}.
Then, there exists an integer $j_0 \geq 1$ (as in Proposition \ref{lemma:ptnotempty}) such that $g_n (O) \subset O$ for all integer $n \geq j_0$.

It is easy to see that
$$\{\lambda \in \mathbb{C},\; \lim g_n (\lambda)=0\}= \bigcup_{n=j_0}^{+\infty}g_n^{-1}(O).$$

Let $z_0$ be a critical points of $g_n,\; n \geq j_0$. By (\ref{ps2}), $g_n= h_n \circ \tilde{f}_{n-1}$, then $\tilde{f}_{k}(z_0)=0$ where $1 \leq k \leq n-1$.
Hence
$ g_{n}(z_0)= h_n \circ f_{n-1} \circ \ldots \circ f_{k+1}(0)$.

Since  $\lim_{j \rightarrow \infty} p_j = 1$ and $p_i > \rho= 2 (\sqrt{2}-1)$ for all $i \geq k+1$, then we have by
the  same argument of proposition \ref{lemma:ptnotempty} that $\lim g_n (z_0)=0$. Hence
 $z_0 \in g_n^{-1}(O)$ for all $n \geq j_0$. Thus we
 deduce by Riemann-Hurwitz formula (see \cite{Mil}), that $g_n^{-1}(O)$ is connected for any integer $n \geq j_0$.
Since $g_n^{-1}(O),\; n \geq j_0$, is a sequence of  increasing sets, we deduce that  $\bigcup_{n=j_0}^{+\infty}g_n^{-1}(O)$ is a connected set. Hence
$$\{\lambda \in \mathbb{C},\; \lim g_n (\lambda)=0\} \subset V.$$

On the other hand,
since  $(g_n)_{n \geq j_0}$ is a  uniformly bounded sequence (by $1$) of holomorphic functions defined on an open subset  $V \subset int(E_{\bar{d},\bar{p}})$. Hence, we deduce by Arzel\`a-Ascoli Theorem (see \cite{Co}), that  $(g_n)_{n \geq j_0}$ is normal in $V $.
That is, there exists a subsequence $(g_{n_k})_{k \geq j_0}$ of $(g_n)_{n \geq j_0}$ such that  $g_{n_k}$ converges to a function $g$ on every compact subset of $V$.

Since $g_n$ converges uniformly on $O$ to $0$, we deduce that $g_n$ converges uniformly on every compact set in $V$ to $g=0$. Hence
$$V \subset \{\lambda \in \mathbb{C},\; \lim g_n (\lambda)=0\}.$$
This ends the proof of Proposition \ref{composant}.
$\square$

\bigskip

The previous result makes the assertion $\sigma_{p} (c_0,\bar{d},\bar{p}) = int (E_{\bar{d},\bar{p}})$ equivalent to $int (E_{\bar{d},\bar{p}})$ is connected. In \cite{msv}, we prove (in collaboration with O. Sester) that if $\bar{d}$ constant and $p_i \geq 2 (\sqrt{2}-1)$ for all $i \geq 1$ then $E_{\bar{d},\bar{p}}$ is a quasidisk. So, in this case $\sigma_{p} (c_0,\bar{d},\bar{p}) = int (E_{\bar{d},\bar{p}})$. We conjecture that this holds for $p_i \geq 2 (\sqrt{2}-1)$ for all $i \geq 1$ even if $\bar{d}$ is non-constant. In this direction, we are only able to show that $E_{\bar{d},\bar{p}}$ is connected, which is the content of our next result:

\bigskip

\begin{proposition}
\label{interi}
Assume that $\lim_{j \rightarrow \infty} p_j = 1$. If $p_i \geq \rho= 2 (\sqrt{2}-1)$ for all $i \geq 1$, then $E_{\bar{d},\bar{p}}$ is connected.
\end{proposition}

\bigskip

\noindent \textbf{Proof:}

By Lemma \ref{cor:bound}, we have
$$ E_{\bar{d},\bar{p}}= \bigcap_{n=1}^{+\infty} g_{n}^{-1} \overline{D(0,1)},$$
where $ g_{n}^{-1} \overline{D(0,1)} \subset  g_{n-1}^{-1} \overline{D(0,1)}$ for all $n \geq 1$.

On the other hand if $R >1$,
it is easy to see that
$$ E_{\bar{d},\bar{p}}= \bigcap_{n=1}^{+\infty} g_{n}^{-1} D(0,R)= \bigcap_{n=1}^{+\infty} g_{n}^{-1} \overline{D(0,R)}, $$
where $ g_{n}^{-1} D(0,R) \subset  g_{n-1}^{-1} D(0,R)$ for all $n \geq 1$.

Let $z_0$ be a critical points of $g_i,\; i \geq 1$, then  we have by
the  same method of proposition \ref{lemma:ptnotempty} that $\lim g_n (z_0)=0$.
Hence $z_0 \in E_{\bar{d},\bar{p}} \subset  g_{i}^{-1} D(0,R)$.
 Thus by  Riemann-Hurwitz formula (see \cite{Mil}, we deduce that  $g_{i}^{-1} D(0,R)$ is connected. Then
 $g_{i}^{-1} \overline{D(0,R)}$ is also connected.
 Hence   $ E_{\bar{d},\bar{p}}$ is the intersection of a decreasing sequence of compact connected sets. Therefore  $ E_{\bar{d},\bar{p}}$ is also connected.
$\square$

\bigskip

\begin{proposition}
In $c$, the point spectrum $\sigma_p (c,\bar{d},\bar{p})$ equals $\sigma_p (c_0,\bar{d},\bar{p}) \cup \{1\}$.
In particular  $\sigma_p (c,\bar{d},\bar{p})=  \{1\}$ if and only if $ (p_n)_{n \geq 0}$ does not converge to $1$.
\end{proposition}

\noindent \textbf{Proof:}
Let $\lambda \in \sigma_p (\Omega,\bar{d})$ and $v= (v_i)_{i \geq 0} \in c$ an eigenvector of $S$ associated to $\lambda$. Then $\lim v_n =l \in \mathbb{C}.$
Consider $N= q_n+ q_{n+2}$ where $n \in \mathbb{N}$, then, $v_{N}= \iota_\lambda(n) \iota_\lambda(n+2)$,
Tending $n$ to infinity, we get $l^2=l$.

{\bf Case 1: $l=1$}. Then $v_{q_1+ q_n}= \iota_\lambda(1) \iota_\lambda(n)$.  Tending $n$ to infinity, we have $\iota_\lambda(1)= 1$,
and thus $\lambda=1$.

{\bf Case 2: $l=0$}, then $\lim p_n=1.$

In both cases, we have $\sigma_p (c,\bar{d},\bar{p})=  \sigma_p (c_0,\bar{d},\bar{p})\cup \{1\} $.
$\square$

\smallskip


\subsection{Case $\Omega= l^1$}

\begin{proposition}
\label{ptspecl1}
If $\prod_{i=1}^{\infty} p_i =0$, then the point spectrum $\sigma_p (l^1,\bar{d},\bar{p}) = \emptyset$.
\end{proposition}

\medskip

\noindent \textbf{Proof:}
Assume that $\prod_{i=1}^{\infty} p_i = 0$ and let $\lambda$ be an eigenvalue of $S_{\bar{d},\bar{p}}$ on $l^1$, then $\sum_{j=0}^{+\infty} \vert \iota_{\lambda}(j+1) \vert $ is convergent.
Since $d_{j} \geq 1$ for all $j \in \mathbb{N}$, we deduce that
$\sum_{j=0}^{+\infty} \vert \iota_{\lambda}(j) \vert^{d_{j}}$ is convergent.

By (\ref{xc}), we deduce that the serie $\sum_{i=0}^{+\infty} 1- p_i $ is convergent, and this contradicts the fact that  $\prod_{i=1}^{\infty} p_i = 0$.
$\square$

\bigskip

As mentioned before, If $\bar{p}$ is constant equal to $1$, then $\sigma_{p} (l^1,\bar{d},\bar{p}) = \mathbb{D}(0,1)$. Thus the condition on $\bar{p}$ in the statement of Proposition \ref{ptspecl1} is necessary. Indeed the condition on the statement is necessary even when $p_j \neq 1$ for every $j \ge 1$. This is shown in the next result.

\medskip

\begin{proposition}
\label{lemma:ptnotemptyl1}
If $\bar{p}$ is monotone increasing and $\prod_{i=1}^{\infty} p_i > 0$, then $ \sigma_{p} (l^1,\bar{d},\bar{p})=\sigma_{p} (c_0,\bar{d},\bar{p}) $. Precisely $ \sigma_{p} (l^1,\bar{d},\bar{p})$ is equal to   the connected component  of  $ int (E_{\bar{d},\bar{p}}) $ that contains $0$.
\end{proposition}

\noindent \textbf{Proof:} The proof follows directly from the next two claims:

\bigskip

\noindent \textbf{Claim 1:} If  $\bar{p}$ is increasing and $\prod_{i=1}^{\infty} p_i > 0$, then $\iota_\lambda \in l^1$.

\bigskip

\noindent \textbf{Claim 2:} 
 $\iota_\lambda \in l^1$, if and only if, $v_\lambda \in l^1$.

\bigskip

\noindent \emph{Proof of Claim 1:} Consider the proof of Proposition \ref{lemma:ptnotempty}. Since $\bar{p}$ is increasing, the choice of $j_k$ implies that
\begin{eqnarray}
\label{pi}
\frac{r^{2\eta^k}}{p_j} + \frac{1-p_{j}}{p_{j}} \ge r^{\eta^{k+1}} \, , \ \textrm{ for } j_{k} \le j \leq  j_{k+1} - 1 \, .
\end{eqnarray}
Thus, for all $k$ sufficiently large,
\begin{equation}
\label{eq:isum}
\sum_{j=j_{k}}^{j_{k+1} -1} \frac{(1-p_j)}{p_{j}} \ge \sum_{j=j_{k}}^{j_{k+1} -1} \frac{r^{2\eta^k}}{p_j} \big( \rho \, r^{(\eta - 2)\eta^k} - 1 \big) \ge \sum_{j=j_{k-1}}^{j_k -1} \frac{r^{2\eta^k}}{p_j}.
\end{equation}
Observe that $\prod_{i=1}^{\infty} p_i > 0$ implies that $\sum_{j\ge 1} (1-p_j) < \infty$ and hence
$\sum_{j\ge 1} \frac{(1-p_{j})}{p_{j}}  < \infty $. Thus   by (\ref{eq:isum}), we have that
$
\sum_{k \ge 0} \sum_{j=j_{k}}^{j_{k+1} -1} \frac{r^{2\eta^k}}{p_j} < \infty .$
Hence $
\sum_{k \ge 0} \sum_{j=j_{k}}^{j_{k+1} -1} \frac{r^{2\eta^k}}{p_j} < \infty .$
Thus, by (\ref{pi}), we have
$$\sum_{k \ge 0} \sum_{j=j_{k}}^{j_{k+1} -1} r^{\eta^{k}} < \infty.$$
By (\ref{eq:iconv}), we obtain that $\iota_\lambda \in l^1$. This finishes the proof of claim 1.

\bigskip

\noindent \emph{Proof of Claim 2:} We prove a general convergence criteria for series: Let $(z_j)_{j \ge 1}$ be a sequence of positive real numbers bounded above by one and define
$$
v(n) = \prod_{j=1}^{\infty} {z_{j}}^{a_{j}(n)} \, , \ n \ge 0.
$$
We have that $(z_j)_{j \ge 1} \in l^1$, if and only if, $(v(n))_{n \ge 0} \in l^1$. Clearly $(z_j)_{j \ge 1} \notin l^1$ implies $(v(n))_{n \ge 0} \notin l^1$. Conversely, assume that  $(z_j)_{j \ge 1} \in l^1$.
We have


\begin{eqnarray}
\label{srl}
\sum_{n \ge 0} |v(n)| = 1 + \sum_{i=1}^{d_{1}-1} z_1^i + \sum_{j \ge 2} \Big( \prod_{k=1}^{j-1}
\sum_{i=0}^{d_{k}-1} z_k^i \Big) \Big( \sum_{i=1}^{d_{j}-1} z_j^i \Big).
\end{eqnarray}
Put $a_j=  \Big( \prod_{k=1}^{j-1}\sum_{i=0}^{d_{k}-1} z_k^i \Big) \Big( \sum_{i=1}^{d_{j}-1} z_j^i \Big)$.

We have $$ \frac{a_{n+1}}{a_n}= B_n \frac{z_{n+1}}{z_{n}},\; \forall n \in \mathbb{N},$$
where
$$
B_n= \frac{(1- z_n^{d_{n}})(1- z_{n+1}^{d_{n+1} -1}) }{(1- z_n^{d_{n}-1})(1- z_{n+1}) }.
$$

Since $\lim z_n = 0$, for $n$ sufficiently large
$$
|B_n| \le  \frac{(1+ |z_n|)(1+ |z_{n+1}|)}{(1- |z_n|)(1- |z_{n+1}|)},
$$
and then, using that $(z_j)_{j \ge 1} \in l^1$, we have that
$$
\prod_{j=1}^{\infty} |B_j| < \infty \, .
$$
Together with the fact that $ a_n= \frac{z_n}{z_1}  B_{n-1} \ldots B_1,\; \forall n \in \mathbb{N},$
we deduce that $(v(n))_{n \ge 0} \in l^1$.
$\square$

\bigskip \bigskip

\medskip

\begin{remark}
\label{rtj}
Using the same proof, we have $\iota_\lambda \in l^{\alpha}, \alpha >1$, if and only if, $v_\lambda \in l^{\alpha}$.
\end{remark}

\vspace{1em}



\bigskip

\begin{theorem}
\label{prop:spectruml1-1}
If $\prod_{i=1}^{\infty} p_i = 0$, then $\sigma_r (l^1,\bar{d},\bar{p})$ contains a countable subset $X$ of the boundary of $E_{\bar{d},\bar{p}}$.
Precisely
$$
X=  \bigcup_{n=1}^{+\infty} \tilde{f}_{n}^{-1} \{1\} \setminus \bigcup_{n=1}^{+\infty} \tilde{f}_{n}^{-1} \{0\} \subset \sigma_r (l^1,\bar{d},\bar{p}).
$$
Moreover, if  $(d_n)_{n \geq 0}$ is bounded and $\lim sup \; p_n <1$, then  $\sigma_r (l^1,\bar{d},\bar{p})= \bigcup_{n=0}^{+\infty} f^{-n} \{1\} \setminus \bigcup_{n=0}^{+\infty} f^{-n} \{0\}$. If $\prod_{i=1}^{\infty} p_i > 0$, then
$$\sigma_r (l^1,\bar{d},\bar{p}) \cap \bigcup_{n=1}^{+\infty} \tilde{f}_{n}^{-1} \{1\} = \emptyset \, . $$
\end{theorem}

\medskip

\begin{remark}
If  $(d_n)_{n \geq 0}$ is bounded and $\lim sup \; p_n <1$, then we can prove (see Proposition \ref{prop:specodeven}) that for a a large class of $(d_{i})_{i \geq 0}$ and $(p_{i})_{i \geq 0}$, we have $\sigma_r (l^1,\bar{d},\bar{p}) = \bigcup_{n=1}^{+\infty} \tilde{f}_{n}^{-1}\{1\}$.
Hence $\sigma_r (l^1,\bar{d},\bar{p})$ is a countable dense subset of the boundary of $E_{\bar{d},\bar{p}}$.
\end{remark}

\bigskip

\noindent \textbf{Proof of Theorem \ref{prop:spectruml1-1}:}
From usual results in operator theory and Proposition \ref{ptspecl1}, we know that
$$
\sigma_r (l^1,\bar{d},\bar{p}) \subset \sigma_p^\prime (l^\oo,\bar{d},\bar{p})
\subset \sigma_r (l^1,\bar{d},\bar{p}) \cup \sigma_{p} (l^1,\bar{d},\bar{p}).$$
Assume that $\prod_{i=1}^{\infty} p_i = 0$, then by Proposition \ref{ptspecl1}, $\sigma_{p} (l^1,\bar{d},\bar{p}) =\emptyset$.
Thus  $\sigma_r (l^1,\bar{d},\bar{p}) = \sigma_p^\prime (l^\oo,\bar{d},\bar{p})$. By Lemmas \ref{lemma:ps} and \ref{lemma:left-eigenvector} and equation (\ref{tb}),
we see that
$$\sigma_p^\prime (l^\oo,\bar{d},\bar{p})=
\left\{ \lambda \in \mathbb{C} \, : (1/ v_{\lambda}(j)))_{j\ge 1}  {\rm {~is e~bounded~}}  {\rm{and~ }}
\iota_{\lambda}(1) =  \sum_{i=1}^{+\infty}
\frac{(1-p_{i+1}) \prod_{j=2}^{i} p_j}{\prod_{r=1}^{i-1} (\iota_{\lambda}(r))^{d_r-1}} \right\},$$
where $(v_{\lambda}(r)))_{r \geq 1}$ is the sequence defined in Lemma \ref {lemma:ps}.

Hence $\sigma_r (l^1,\bar{d},\bar{p}) $ is contained in the set
\begin{eqnarray}
\label{bla}
E_{\bar{d},\bar{p}} \cap \left \{\lambda \in \mathbb{C} \, : \ (1/ \iota_{\lambda}(j))_{j\ge 1} \mbox
{ is  bounded and }
\iota_{\lambda}(1) =  \sum_{i=1}^{+\infty} \frac{(1-p_{i+1}) \prod_{j=2}^{i} p_j} {\prod_{r=1}^{i-1} (\iota_{\lambda}(r))^{d_r-1}} \right\}
\end{eqnarray}
On the other hand, by (\ref{xc}) and since
$ i_{\lambda}(n)= \frac{1}{p_{n}} \tilde{f}_{n-1}(\lambda)- \frac{1- p_{n}}{p_{n}}$ for all integer $n \geq 1$, we deduce that
\begin{eqnarray}
\label{rth}
\tilde{f}_{n-1}(\lambda)=
i_{\lambda}(n-1)^{d_{n-1}},\; \forall n \geq 1.
\end{eqnarray}

Let $n \in \mathbb{N}$ and $ E_n= \{\lambda \in \mathbb{C},\; \iota_{\lambda}(n)=1\}$.

By (\ref{rth}), we have
\begin{eqnarray}
\label{egal}
  \bigcup_{n=1}^{+\infty} E_n= \bigcup_{n=1}^{+\infty} \tilde{f}_{n}^{-1}\{1\}.
  \end{eqnarray}

Now
assume that there exists  $n_0 \in \mathbb{N}$ and $\lambda \in E_{n_0}$. Then
\begin{eqnarray}
\label{chch}
 \iota_{\lambda}(k)= 1, \; \forall k \geq n_0.
 \end{eqnarray}

 \vspace{0.5em}

Assume that  $\lambda \in \bigcup_{n=1}^{+\infty} \tilde{f}_{n}^{-1} \{1\} \setminus \bigcup_{n=1}^{+\infty}\tilde{f}_{n}^{-1} \{0\}$, then
$\iota_{\lambda}(i) \ne 0$ for all $i <n_0$. From  (\ref{chch}) and (\ref{egal}), we have that $(\iota_{\lambda}(n))_{n \geq 0} \mbox { and } (1/ \iota_{\lambda}(n))_{n \geq 0} \mbox { are bounded  }.$
Moreover,  by (\ref{xc})
 we have

\begin{eqnarray*}
\iota_{\lambda}(1) =   \sum_{i=1}^{+\infty} \frac{(1-p_{i+1}) \prod_{j=2}^{i} p_j}  { \prod_{r=1}^{i-1} (\iota_{\lambda}(r))^{d_r-1} }
&\Longleftrightarrow & \iota_{\lambda}(2)  =  \sum_{i=2}^{+\infty} \frac{(1-p_{i+1}) \prod_{j=3}^{i} p_j}{ \prod_{r=2}^{i-1} (\iota_{\lambda}(r))^{d_{r-1}}}\\
 &\Longleftrightarrow &  \iota_{\lambda}(n_{0}) =  \sum_{i=n_{0}}^{+\infty} \frac{ (1-p_{i+1}) \prod_{j=n_{0}+1}^{i} p_j} { \prod_{r=n_{0}}^{i-1} (\iota_{\lambda}(r))^{d_{r-1}}}
 \\
 &\Longleftrightarrow & 1= \sum_{i=n_0}^{+\infty}  (1-p_{i+1}) \prod_{j=n_{0}+1}^{i} p_j.
 \end{eqnarray*}
 Hence

 \begin{eqnarray}
 \label{zer}
\iota_{\lambda}(1)  =   \sum_{i=1}^{+\infty} \frac{(1-p_{i+1}) \prod_{j=2}^{i} p_j}  { \prod_{r=1}^{i-1} (\iota_{\lambda}(r))^{d_r-1} }
&\Longleftrightarrow  & 0=  \prod_{i=n_{0}+1}^{+\infty} p_i.
 \end{eqnarray}
From this $\lambda \in \sigma_r (l^1,\bar{d},\bar{p})$. Hence $\bigcup_{n=1}^{+\infty} \tilde{f}_{n}^{-1}\{1\}  \setminus \bigcup_{n=1}^{+\infty} \tilde{f}_{n}^{-1} \{0\} \subset \sigma_r (l^1,\bar{d},\bar{p})$.

By (\ref{zer}), we deduce that if $\prod_{i=1}^{\infty} p_i > 0$, then $\sigma_r (l^1,\bar{d},\bar{p}) \bigcap \bigcup_{n=0}^{+\infty} \tilde{f}_{n}^{-1}\{1\}  = \emptyset$.

\bigskip

It remains to prove the following result.

\begin{proposition}
\label{prop:spectruml1-2}
If  $(d_n)_{n \geq 0}$ is bounded and $\lim sup \; p_n <1$, then $  \sigma_r (l^1,\bar{d},\bar{p}) \subset \bigcup_{n=0}^{+\infty} \tilde{f}_{n}^{-1} \{1\}\setminus \bigcup_{n=1}^{+\infty} \tilde{f}_{n}^{-1} \{0\}$.
\end{proposition}

\medskip


\bigskip

\noindent \textbf{Proof:}
Let $\lambda \in E_{\bar{d},\bar{p}}$ and $\lambda \not \in \bigcup_{r=1}^\infty \tilde{f}_r^{-1} \{1\}$.
Then $ \vert i_\lambda (n)\vert \leq 1$ and  $  i_\lambda (n)\ne  1$ for all integer $n \geq 1$.

{\bf Case 1: $\liminf_{j \ra \oo} |\iota_\lambda (j)|= C < 1$.}

Let $\varepsilon >0$ such that $C+\varepsilon <1$.  Then there exists an increasing  sequence of positive integers $(k_{j})_{j \geq 1}$ such that $ |\iota_\lambda (k_j)|\leq C+\varepsilon$ for all $j \geq 1$.

Now, consider the sequence $(x_n)_{n \geq 1}$ defined by
$x_n = q_{k_{1}}+ \cdots + q_{k_{n}}$.
Hence $\vert v_\lambda (x_n)\vert = \prod_{r=1}^{n} \vert (\iota_{\lambda}(k_{r}))\vert \leq (C+\varepsilon)^n$.
Thus $\vert v_\lambda (x_n)\vert $ converges to $0$ as $n$ goes to infinity, and $\frac{1}{\vert v_\lambda (x_n)\vert} $ is not bounded.
Thus $\lambda \not \in \sigma_r (l^1,\bar{d},\bar{p})$.

{\bf Case 2: $\liminf_{j \ra \oo} |\iota_\lambda (j)|= 1$}.
Then $\lim_{n \ra \oo} |\iota_\lambda (n)|= 1$.

For all integer $n \geq 1$, put $\iota_\lambda (n)= r_n e^{i \theta_n}$ where $0 \leq r_n \leq 1$ and $\theta_n \in [0, 2 \pi).$

By (\ref{xc}), we get for all integer $n \geq 1$,

\begin{eqnarray}
\label{cose}
r_{n+1} \cos \theta_{n+1} = \frac {r_{n}^{d_{n}} \cos d_n \theta_{n}}{p_{n+1}} - \frac{1- p_{n+1}}{p_{n+1}},\;
r_{n+1} \sin \theta_{n+1} = \frac{r_{n}^{d_{n}} \sin d_n \theta_{n}}{p_{n+1}}
 \end{eqnarray}
 Thus
 \begin{eqnarray}
 \label{vbx}
  p_{n+1}^2 r_{n+1}^2=  r_{n}^{2 d_{n}} + (1- p_{n+1})^2 -  2 (1- p_{n+1})r_{n}^{d_{n}} \cos d_n \theta_{n}.
 \end{eqnarray}

Let $\varepsilon >0$, then by (\ref{vbx}) and the fact that $\lim  r_n= 1$ and $(d_n)$ bounded, we deduce that
  there exists an integer $N$ such that for all integer $n \geq N$, we have
$2 (1- p_n) (1- \cos d_n \theta_{n}) < \varepsilon.$

Since $\lim sup \; p_n < 1$ for all $n$, then $\cos d_n \theta_{n}$ converges to $1$.
Hence, by (\ref{rth}), $\lim \tilde{f}_{n-1}(\lambda) = \lim \iota_\lambda (n-1)^{d_{n-1}} =1.$
Thus $\lim \iota_\lambda (n) =1.$
Then, given $0 <\varepsilon <1$, there exists an integer $n_0 \geq 1$ such that for all integer
 $n \geq n_0$, we have
$$\vert 1+ \iota_\lambda (n)+ \cdots \iota_\lambda (n)^{d_{n}-1} \vert \geq d_n - (d_n -1) \varepsilon \geq 2- \varepsilon.$$
Therefore for all integer $n \geq n_0$, we have
 $ \vert \iota_{\lambda} (n+1) -1 \vert = \frac{1}{p_{n+1}} \vert  \iota_\lambda (n)^{d_{n}}-1 \vert \geq
 \frac{ 2- \varepsilon}{p_{n+1}}  \vert \iota_{\lambda} (n) -1 \vert \geq \frac{(2-  \varepsilon)^{n-n_0+1}}{p_{n+1} \ldots p_{n_0}}  \vert \iota_{\lambda} (n_{0}) -1 \vert$.
 Hence $\iota_{\lambda} (n_{0}) =1$. Thus $\tilde{f}_{n_0}(\lambda)=1,$
 which is a contradiction.
 This ends the proof of Proposition \ref{prop:spectruml1-2} and Theorem \ref{prop:spectruml1-1}. $\square$


\bigskip \medskip

\begin{conjecture}
\label{residspect}
In $l^1,\;
\sigma_r (l^1,\bar{d},\bar{p}) = \emptyset$ if $\prod_{i=1}^{\infty} p_i > 0$, and $ \sigma_r (l^1,\bar{d},\bar{p})=  \bigcup_{n=1}^{+\infty} \tilde{f}_{n}^{-1}\{1\}  \setminus \bigcup_{n=1}^{+\infty} \tilde{f}_{n}^{-1} \{0\}$ otherwise.
\end{conjecture}

\medskip

\begin{proposition}
\label{prop:specodeven}
 Assume that   $(d_n)_{n \geq 0}$ is bounded and $\lim sup \; p_n <1$. Then  the following properties are valid.
 \begin{enumerate}
 \item
 If all $d_k,\; k \geq 1$ are odd, then $\sigma_r (l^1,\bar{d},\bar{p})= \bigcup_{n=0}^{+\infty} \tilde{f}_{n}^{-1}\{1\}  $.
 \item
 If all $k \ge 1$, $d_k$ is even and $p_k > \frac{1}{2}$,  then $\sigma_r (l^1,\bar{d},\bar{p})= \bigcup_{n=0}^{+\infty} \tilde{f}_{n}^{-1}\{1\}$.
 \item
 If there exists an integer $k \geq 1$ such that $d_{k}$ is even, then
 \begin{enumerate}
 \item
 If $p_{k}= \frac{1}{2}$, then $\sigma_r (l^1,\bar{d},\bar{p}) \ne \bigcup_{n=0}^{+\infty} \tilde{f}_{n}^{-1}\{1\}  $.
\item
 If $p_{k}< \frac{1}{2}$ and $d_{k-1}$ is even, then there exists $a \in (0,1)$ such that $p_{k-1}=a$ implies that $\sigma_r (l^1,\bar{d},\bar{p}) \ne \bigcup_{n=0}^{+\infty} \tilde{f}_{n}^{-1}\{1\}$.
 \end{enumerate}
 \item Consider $\bar{p}= (p_i)_{i \geq 0}$ random such that $p_i$'s are iid random variables with continuous distribution. Then, given any sequence $\bar{d}=(d_k)_{k \geq 1}$ of integers such that $d_0=1$ and $d_k \geq 2$ for all $k$, we have $\sigma_r (l^1,\bar{d},\bar{p})= \bigcup_{n=0}^{+\infty} \tilde{f}_{n}^{-1}\{1\}$ with probability one.
\end{enumerate}
\end{proposition}

\bigskip

\noindent \textbf{Proof:} First note that (2) is a direct consequence of the following simple assertion:  If $d_k,\; k \geq 0$ is even and $p_k > \frac{1}{2}$, then 
$\iota_\lambda(k) \in (-1,1)$ implies that $\iota_\lambda(k+1) \in (-1,1)$. In particular, if $d_k$ is even and $p_k > \frac{1}{2}$, for all $k$, and $\iota_\lambda(k_0) = 0$ for some $k_0$, then $\iota_\lambda(k) \in (-1,1)$ for all $k \ge k_0$. This means that
$\bigcup_{n=0}^{+\infty} \tilde{f}_{n}^{-1}\{1\}  \cap \bigcup_{n=1}^{+\infty} \tilde{f}_{n}^{-1} \{0\} = \emptyset$.

 By Theorem \ref{prop:spectruml1-1}, $\sigma_r (l^1,\bar{d},\bar{p})= \bigcup_{n=0}^{+\infty} \tilde{f}_{n}^{-1}\{1\} \setminus \bigcup_{n=0}^{+\infty} \tilde{f}_{n}^{-1}\{0\} $.
Let us analyze when $\bigcup_{n=0}^{+\infty} \tilde{f}_{n}^{-1}\{1\} \cap \bigcup_{n=0}^{+\infty} \tilde{f}_{n}^{-1}\{0\}$ is empty or not. Let    $\lambda \in \bigcup_{n=0}^{+\infty} \tilde{f}_{n}^{-1}\{1\} \cap \bigcup_{n=0}^{+\infty} \tilde{f}_{n}^{-1}\{0\}$, then there exist integers $1 \leq m < n$ such that  $\iota_{\lambda}(m) =0$ and  $\iota_{\lambda}(n) =1$. Thus there exists an integer $k$ such that $i_{\lambda}(k) \ne  1$ and $i_{\lambda}(k+1)= 1$.
By (\ref{xc}), we obtain $i_{\lambda}(k)^{d_{k}} =  1$.

Assume that $d_{k}$ is odd,  then by (\ref{xc}), $i_{\lambda}(j) \not \in \mathbb{R}$  for all $0 \leq j \leq k$, absurd. Then if all  $d_i,\; i \geq 0$ are odd, then $\bigcup_{n=0}^{+\infty} \tilde{f}_{n}^{-1}\{1\}\cap \bigcup_{n=0}^{+\infty} \tilde{f}_{n}^{-1}\{0\}= \emptyset $ and we obtain (1).

Now, suppose that  $d_{k}$  is even. Choose $\lambda$ such that $\iota_\lambda(k) = - 1$ which implies $i_{\lambda}(k)^{d_{k}}= 1$, i.e. $\lambda \in \tilde{f}_{k}^{-1}\{1\}$. Thus $i_{\lambda}(k-1)^{d_{k-1}}= p_k i_{\lambda}(k) - (1-p_k) = 1- 2 p_{k}$.

If $ p_{k}= 1/2$, then $i_{\lambda}(k-1)= 0$. Hence  $\sigma_r (l^1,\bar{d},\bar{p}) \cap \bigcup_{n=0}^{+\infty} \tilde{f}_{n}^{-1}\{1\} \ne \emptyset $ and we obtain (3.a).

If $p_k < 1/2$ and $d_{k-1}$ is even, then we can restrict the choice of $\lambda$ so that $i_{\lambda}(k-1) = \sqrt[d_{k-1}]{1- 2 p_{k}} \in (-1,0)$. Now take $p_{k-1}$ such that $i_{\lambda}(k-1)= -\frac{1- p_{k-1}}{p_{k-1}}$ and we get  $i_{\lambda}(k-2)=0$. Hence we deduce (3.b).

It remains to prove (4). It is enough to show that
\begin{equation}
\label{eq:(4)}
P\Big( \bar{p} \, : \, \tilde{f}_{k}^{-1}\{1\} \cap \bigcup_{n=1}^{k-1} \tilde{f}_{n}^{-1}\{0\} \neq \emptyset \Big) = 0 \, ,
\end{equation}
for every $k \ge 2$. This holds because given $p_k$ there is only finite possible choices of $p_1$, ..., $p_{k-1}$, that implies $\tilde{f}_{k}^{-1}\{1\} \cap \bigcup_{n=1}^{k-1} \tilde{f}_{n}^{-1}\{0\} \neq \emptyset$. Since the random vector $(p_1,...,p_{k-1})$ has continuous distribution and is independent of $p_k$, it will take values in a finite set with probability zero. Therefore we have \ref{eq:(4)}. $\square$

\medskip

\begin{remark}
If $d_k= 2$ for all $k$ and $p_k = a$, then $\sigma_r (l^1,\bar{d},\bar{p})= \bigcup_{n=1}^{+\infty} \tilde{f}_{n}^{-1}\{1\} $ if $a \ne \frac{1}{2} $ and  $\sigma_r (l^1,\bar{d},\bar{p})= \{1\}$ otherwise. This  last case corresponds to the case where the Julia set $ E_{\bar{d},\bar{p}}$ is a dendrite.
It will be interesting to characterize in the general case, the relation between the fact that  $\sigma_r (l^1,\bar{d},\bar{p}) \ne \bigcup_{n=0}^{+\infty} \tilde{f}_{n}^{-1}\{1\}  $ and topological properties of $ E_{\bar{d},\bar{p}}$.
\end{remark}

\subsection{Case $\Omega= l^{\alpha},\; \alpha >1 $}

\vspace{2em}

In this section we consider $\alpha > 1$ fixed. We proceed in analogy to the case $l^1$, obtaining versions of Propositions \ref{ptspecl1} and \ref{lemma:ptnotemptyl1}.

\medskip

\begin{proposition}
If $\sum_{j=1}^\infty (1-p_j)^\alpha = \infty$, then $\sigma_p (l^\alpha,\bar{d},\bar{p})= \emptyset$.
\end{proposition}
\noindent \textbf{Proof:}
Take $\lambda \in E_{\bar{d},\bar{p}}$. We have that $\lambda \in \sigma_p (l^\alpha,\bar{d},\bar{p})$, if and only if, $v_\lambda \in l^{\alpha}$, which, by remark \ref{rtj}, is equivalent to $\iota_\lambda \in l^{\alpha}$. Then, if we suppose that $\lambda \in \sigma_p (l^\alpha,\bar{d},\bar{p})$, we have that $(\iota_\lambda(j)^{d_j})_{j\ge 1}$ and $(p_j \, \iota_\lambda(j+1))_{j \ge 1}$ are in $l^\alpha$. Since $(1-p_{j+1}) = \iota_\lambda(j)^{d_j} -  p_{j+1} \iota_\lambda(j+1)$, we have that $(1-p_j)_{j\ge 1}\in l^{\alpha}$. Therefore $(1-p_j)_{j\ge 1}\notin l^{\alpha}$ implies that $\sigma_p (l^\alpha,\bar{d},\bar{p}) = \emptyset$.
$\square$

\bigskip

\begin{proposition}
If $\bar{p}$ is monotone increasing and $\sum_{j=1}^\infty (1-p_j)^\alpha < \infty$, then $\sigma_p (l^\alpha,\bar{d},\bar{p})$ equals to the connected component of $int (E_{\bar{d},\bar{p}}) $ that contains $0$.
\end{proposition}
\noindent \textbf{Proof:}
Consider the proof of Proposition \ref{lemma:ptnotempty}. Since $\bar{p}$ is increasing, the choice of $j_k$ implies that
\begin{eqnarray}
\label{xcj}
2^\alpha \Big( \frac{r^{2 \alpha \eta^k}}{p^\alpha_j} + \Big( \frac{1-p_{j}}{p_{j}} \Big)^\alpha \Big) \ge r^{\alpha \eta^{k+1}} \, , \ \textrm{ for } j_{k} \le j < j_{k+1} - 1 \, .
\end{eqnarray}

Thus, for all $k$ sufficiently large,
\begin{eqnarray*}
\sum_{j=j_{k}}^{j_{k+1} -1} \Big( \frac{1-p_j}{p_{j}} \Big)^\alpha & \ge &
\sum_{j=j_{k}}^{j_{k+1} -1} 2^{-\alpha} \, \frac{r^{2 \alpha \eta^k}}{p_j^{\alpha}} \big( \rho^{\alpha} \; r^{\alpha(\eta - 2)\eta^k} - 2^ \alpha \big) \\
& \ge & 2^{-\alpha} \sum_{j=j_{k}}^{j_{k+1} -1} \frac{r^{2 \alpha \eta^k}}{p_j^ {\alpha}}.
\end{eqnarray*}
Since $\sum_{j\ge 1} (1-p_j)^\alpha < \infty$ implies $\sum_{j\ge 1} \Big( \frac{1-p_j}{p_{j}} \Big)^\alpha  < \infty $, we have that
$\sum_{k \ge 0} \sum_{j=j_{k}}^{j_{k+1} -1} \frac{r^{2 \alpha \eta^k}}{p_j^ {\alpha}} < \infty .$
Hence by (\ref{xcj}), we have
$$\sum_{k \ge 0} \sum_{j=j_{k}}^{j_{k+1} -1} r^{\alpha \eta^{k}} < \infty.$$
By (\ref{eq:iconv}), we obtain that $\iota_\lambda \in l^\alpha$.
$\square$

\bigskip \medskip

\section{Proofs of the Technical Lemmas}
\label{sec:proofs}

\bigskip

\noindent \textbf{Proof of Lemma \ref{bound}:}
Take $p_j\in (0,1)$ and $z\in \mathbb{C}$ with $|z|>1$ then
\begin{equation}
\label{lowerbound1}
\left| \frac{z-(1-p_j)}{p_j}\right| \ge \frac{|z| - (1-p_j)}{p_j} = \frac{|z|-1}{p_j} + 1 > |z| > 1 \, .
\end{equation}
Thus, we obtain, for every $z\in \mathbb{C}$ with $|z|>1$ and $j\ge 1$, that
$$
|f_{j}(z)| > |z|^{d^j} \, .
$$

Now suppose $|\tilde{f}_{r}(z)|>1$ for some $r>1$, then by induction one can show that for $j>r$
\begin{equation}
\label{lowerbound2}
|\tilde{f}_{j}(z)| \ge |\tilde{f}_{r}(z)|^{d_j ... d_{r+1}} \, .
\end{equation}
Indeed
$$
|\tilde{f}_{j+1}(z)| = \left| \frac{\tilde{f}_{j}(z)-(1-p_{j+1})}{p_{j+1}} \right|^{d_{j+1}} \ge \left| \frac{|\tilde{f}_{j}(z)|-1}{p_{j+1}} + 1 \right|^{d_{j+1}} \ge |\tilde{f}_{j}(z)|^{d_{j+1}} \, .
$$
From (\ref{lowerbound2}) we see that $\lim_{j\ra +\oo} |\tilde{f}_{j}(z)| = +\oo$ whenever $|\tilde{f}_{r}(z)|>1$ for some $r>1$.
In particular, if $|\tilde{f}_{r}(z)|>1$ then $z \notin E_{\bar{d},\bar{p}}$.\\
Now, suppose $|z-(1-p_1)|>p_1,$ this implies that $|f_1(z)|>1$ and then $z\notin E_{\bar{d},\bar{p}}$. Analogously, if $|\tilde{f}_j(z) - (1-p_{j+1})| > p_{j+1}$, we have that $|\tilde{f}_{j+1}(z)|>1$ and then $z \notin E_{\bar{d},\bar{p}}$.
$\square$

\bigskip \bigskip

\noindent \textbf{Proof of Lemma \ref{lemma:ps}:}
Let $v=(v_n)_{n\ge 0}$ be a sequence of complex numbers and suppose that $(S \, v)_n = \lambda \, v_n$ for every $n \ge 0$. We shall prove that $v$ satisfies (\ref{ps1}). The proof is based on the following representation
\begin{eqnarray}
\label{ps4}
(S \, v)_n & = & \left( \prod_{j=1}^{\zeta_{n}} p_j \right) \, v_{n+1} + (1-p_1) v_n \nn \\
& & + \sum_{r=1}^{\zeta_{n}-1}
\left( \prod_{j=1}^{r} p_j \right) (1-p_{r+1}) v_{n-\sum_{j=1}^r (d_j-1) q_{j-1}} \, ,
\end{eqnarray}
for $\zeta_{n} \geq 2$
and $(S \, v)_n = p_1 v_{n+1} + (1-p_1) v_n$ if $\zeta_{n}= 1$. This representation follows directly from
the definition of the transition probabilities in (\ref{tp1}). From (\ref{ps4}), we show (\ref{ps1}) by induction.

Indeed, for $n=1$ we have that
$$
\lambda v_0 = (1-p_1)v_0 + p_1 v_1 \ \Rightarrow \ v_1 = \left( \frac{\lambda - (1-p_1)}{p_1}\right) v_0 = \iota_{\lambda}(1) \, v_0 \, .
$$
Now fix $n \ge 1$ and suppose that (\ref{ps1}) holds for every $1 \le j \le n$. By (\ref{ps4}), since $(S \, v)_n = \lambda \, v_n$, we have that
\begin{equation}
\label{ps5}
\frac{v_{n+1}}{v_0 \, \prod_{r=\zeta_{n}+1}^{\infty} (\iota_{\lambda}(r))^{a_{r}(n)}}
\end{equation}
is equal to
\begin{eqnarray}
\label{ps6}
\lefteqn{\!\!\!\!\!\!\!\!\!
\frac{ [\lambda - (1-p_1)] \Big[ \prod_{r=1}^{\zeta_{n}-1} (\iota_{\lambda}(r))^{d_r-1} \Big] (\iota_{\lambda}(\zeta_{n}))^{a_{\zeta_{n}}(n)} }
{ \prod_{j=1}^{\zeta_{n}} p_j }
} \nn \\
& & - \frac{ (1-p_2) \Big[ \prod_{r=2}^{\zeta_{n}-1} (q_{\lambda}(r))^{d_r-1} \Big] (\iota_{\lambda}(\zeta_{n}))^{a_{\zeta_{n}}(n)} }
{ \prod_{j=2}^{\zeta_{n}} p_j }
\nn \\
& & ... \, - \frac{1-p_{\zeta_{n}}}{p_{\zeta_{n}}} (\iota_{\lambda}(\zeta_{n}))^{a_{\zeta_{n}}(n)} \, .
\end{eqnarray}
Since
$$
\iota_{\lambda}(1) = \frac{\lambda - (1-p_1)}{p_1} \, ,
$$
the first term in (\ref{ps6}) is equal to
$$
\frac{ \iota_{\lambda}(1)^{d_1} \, \Big[ \prod_{r=2}^{\zeta_{n}-1} (\iota_{\lambda}(r))^{d_r-1} \Big] (\iota_{\lambda}(\zeta_{n}))^{a_{\zeta_{n}}(n)} }
{ \prod_{j=2}^{\zeta_{n}} p_j } \, .
$$
Summing with the second term we get
$$
\left( \frac{ (\iota_{\lambda}(1))^{d_1} - (1-p_2)}{p_2} \right)
\frac{ \Big[ \prod_{r=2}^{\zeta_{n}-1} (\iota_{\lambda}(r))^{d_r-1} \Big] (\iota_{\lambda}(\zeta_{n}))^{a_{\zeta_{n}}(n)} }
{ \prod_{j=3}^{\zeta_{n}} p_j } \, ,
$$
which is equal to
$$
\frac{ (\iota_{\lambda}(2))^{d_2} \, \Big[ \prod_{r=3}^{\zeta_{n}-1} (\iota_{\lambda}(r))^{d_r-1} \Big] (\iota_{\lambda}(\zeta_{n}))^{a_{\zeta_{n}}(n)} }
{ \prod_{j=3}^{\zeta_{n}} p_j } \, .
$$
By induction we have that the sum of the first $\zeta_{n}-1$ terms in (\ref{ps6}) is equal to
$$
\frac{ (\iota_{\lambda}(\zeta_{n}-1))^{d_{\zeta_n - 1}} \, (\iota_{\lambda}(\zeta_{n}))^{a_{\zeta_{n}}(n)} }
{ p_{\zeta_{n}} } \, .
$$
Finally, summing the previous expression with the last term in (\ref{ps6}) we have that (\ref{ps5}) is equal to
$$
\frac{ (\iota_{\lambda}(\zeta_{n}-1))^{d_{\zeta_n - 1}} - (1-p_{\zeta_{n}}) }{ p_{\zeta_{n}} } \,
(\iota_{\lambda}(\zeta_{n}))^{a_{\zeta_{n}}(n)} = (\iota_{\lambda}(\zeta_{n}))^{a_{\zeta_{n}}(n)+1} \, ,
$$
Therefore,
\begin{eqnarray*}
v_{n+1} & = & v_0 \, (\iota_{\lambda}(\zeta_{n}))^{a_{\zeta_{n}}(n)+1} \, \prod_{r=\zeta_{n}+1}^{\infty} (\iota_{\lambda}(r))^{a_{r}(n)} \\
& = & \prod_{r=1}^{\infty} (\iota_{\lambda}(r))^{a_{r}(n+1)}  \, ,
\end{eqnarray*}
which, by induction, completes the proof of Claim 2.
$\square$

\bigskip \bigskip

\noindent \textbf{Proof of Lemma \ref{lemma:spectrumequaljulia}:}
Consider $\lambda \in E_{\bar{d},\bar{p}}$. Assume that $\lambda \not \in \sigma_p(l^\alpha,\bar{d},\bar{p})$. We will prove that $\lambda$ belongs to the approximate point spectrum of $S_{\bar{d},\bar{p}}$.
  For all  integers $k \geq 2,$ put $w^{(k)}= (v_{\lambda}(0),v_{\lambda}(1), \ldots, v_{\lambda}(k), 0 \ldots 0, \ldots)^{t} \in l^{\alpha} $
where $(v_{\lambda}(r)))_{r \geq 1}$ is the sequence defined in Lemma \ref {lemma:ps}.
Let $u^{(k)}= \frac{w^{(k)}} {\vert \vert  w^{(k)} \vert \vert_{\alpha}}$, then we have the following claim.

\vspace{0.5em}

\noindent {\bf Claim:} $ \lim_{n\rightarrow+\infty} \vert \vert (S- \lambda I) u^{(q_n)}\vert \vert_{\alpha}=0$ where $ q_n= d_0 \ldots d_n$.\\

We assume that $\alpha >1$ (the case $\alpha=1$ can be done using the same method).

\noindent  Indeed, we have
$$\forall i \in \{0,\ldots,k-1\},  ~~\left ((S- \lambda I) u^{(k)}\right)_i=0.$$
 Thus
\begin{eqnarray*}
\sum_{i=0}^{+\infty} \left \vert {((S- \lambda I) u^{(k)})}_i \right \vert ^{\alpha}
 =
\frac { \sum_{i=k}^{+\infty} \left
\vert  \sum_{j=0}^{k}( S - \lambda I)_{i,j}  w^{(k)}_{j} \right \vert }{ \vert \vert w^{(k)}\vert \vert_{\alpha}^{\alpha}}^{\alpha}.
\end{eqnarray*}

Put $a_{i,j}= \vert (S -\lambda I)_{i,j}\vert$ for all $i,j$.
 Let  $\alpha'$ be a conjugate of $\alpha$, i.e, $ \frac1{\alpha}+\frac1{\alpha'}=1.$ Then,
by H\"older inequality we get
$$
\sum_{j=0}^{+\infty}{a}_{i,j} \vert w^{(k)}_{j} \vert  =
\sum_{j=0}^{+\infty} a_{i,j}^{\frac{1}{\alpha'}}  (a_{i,j}^{\frac{1}{\alpha}} \vert w^{(k)}_{j}
\vert )
\leq {\left(\sum_{j=0}^{+\infty}{{a}_{i,j}}\right)}^{\frac{1}{\alpha'}}
\left(\sum_{j=0}^{+\infty}a_{i,j}\vert {w^{(k)}_{j}}\vert ^{\alpha}  \right)^{\frac{1}{\alpha}}
$$



Thus

\begin{eqnarray*}
\left \vert  \sum_{j=0}^{k}( S - \lambda I)_{i,j}  w^{(k)}_{j} \right \vert  ^{\alpha}   \leq C \sum_{j=0}^{k} \vert (S - \lambda I)_{i,j}\vert  \vert w^{(k)}_{j} \vert^{\alpha}
\end{eqnarray*}

where
$C=  \sup_{i \in \N} \left(\sum_{j=0}^{\infty} \vert (S- \lambda I)_{i,j} \vert
 \right)^{\frac {\alpha} {\alpha'}}$ and $\alpha'$
is the conjugate of  $\alpha .$

Observe that $C$ is a finite non-negative constant because $S$ is a stochastic matrix and $\lambda$ belongs to $E$ which is a bounded set.

In this way we have

\begin{eqnarray*}
 \left \vert \left \vert {(S_p- \lambda I) u^{(k)}} \right \vert \right \vert^{{\alpha}}_{\alpha}
& \leq &
C  \sum_{i=k}^{+\infty}\frac{\left( \sum_{j=0}^{k}  \vert w^{(k)}_j \vert^{\alpha}  \vert (S_p- \lambda I)_{ij}\vert  \right)}
{\vert \vert  w^{(k)}\vert \vert_{\alpha}^{\alpha}}\\
& = &
 \frac{C}{\vert \vert  w^{(k)}\vert \vert_{\alpha}^{\alpha}}  \sum_{j=0}^{k} \vert w^{(k)}_j \vert^{\alpha}
\sum_{i=k}^{+\infty}    \vert (S_p- \lambda I)_{ij}\vert.
\end{eqnarray*}

\noindent Now, for $k=q_n$, we will compute the following terms $$ A_{kj}= \sum_{i=k}^{+\infty}  \vert (S_p- \lambda I)_{ij} \vert,\; 0 \leq j \leq k.$$

Assume that $0 \leq j <k =q_n.$ Then $\left(S_p- \lambda I\right)_{ij}=( S_p)_{ij}$ for all $i \geq k$.

{\bf Case 1}:  $j=r \mod d_1,\; 0 < r <d_1$ . Then by (\ref{tp1}), $(S_p)_{ij} \ne 0$ if and only $i=j-1$ or $i=j$.
Hence $(S_p)_{ij}=0$ for all $i \geq k$. Thus
\begin{eqnarray}
\label{t0}
A_{kj}=0.
\end{eqnarray}

{\bf Case 2}:   $j=0$ . Then by  ( \ref{tp1}), we have

\begin{eqnarray}
\label{tts}
A_{kj}= A_{q_n, 0}= \sum_{i=q_n}^{+\infty}   (S_p)_{i0}=  \sum_{i=n+1}^{+\infty} (1-p_{i+1}) \prod_{j=1}^{i} p_j.
\end{eqnarray}
Observe that $ \lim A_{q_n, 0}= 0$.

{\bf Case 3}:   $ j=0 \mod d_1 $ is even and $j >0$. Then $j= a _{n-1}\ldots
a _{s}\underbrace{0 \ldots 0}_{s}=  \sum _{i=s}^{n-1} a_i q_i$  with $s \geq 1$
and  $a_{s}>0$. But by  (\ref{tp1}), $(S_p)_{ij} \ne 0$ if and only if $i=
 a _{n-1}\ldots
a _{s}\underbrace{0 \ldots 0}_{s-m+1}\underbrace{(d-1) \ldots (d-1)}_{m-1}= q_m -1+j$ where $1 \leq m  \leq s$. Hence $i <q_n= k.$

Therefore, in this case
\begin{eqnarray}
\label{tr}
A_{kj}=0.
\end{eqnarray}

Now assume $j=k=q_n$. In this case, we have
$A_{kj}=\vert 1-p-\lambda\vert+  \sum_{i=q_n+1}^{+\infty}  (S_p)_{i, q_n}.$
On the other hand, by (\ref{tp1}), we deduce that
$(S_p)_{i, q_n} \ne 0$ if and only if $i= q_n+ q_m-1$ where $0 \leq m \leq n$
and $(S_p)_{q_n+ q_m-1, q_n} = (1-p_{m+1}) \prod_{j=1}^{m} p_j.$
Therefore
$$
A_{kj}= \sum_{i=q_{n}}^{+\infty}  \vert (S- \lambda I)_{i,q_{n}}\vert= \vert 1-p-\lambda \vert+ \sum_{m=0}^{n}
 (1-p_{m+1}) \prod_{j=1}^{m} p_j.
$$
Hence
\begin{eqnarray}
\label{tt2}
A_{q_n, q_n}=  \vert 1-p-\lambda \vert+ 1-   \prod_{j=1}^{n} p_j.
\end{eqnarray}

By (\ref{t0}),(\ref{tts}),(\ref{tr}) and (\ref{tt2}), we have for $k=q_n$ and $0 \leq j \leq k$,
 \begin{eqnarray}
\label{fro}
A_{kj} \ne 0 \Longleftrightarrow j=0 \mbox { or } j=k=q_n.
\end{eqnarray}

Consequently

\begin{eqnarray*}
\left \vert \left \vert {(S- \lambda I) u^{(q_n)}} \right \vert \right \vert^{\alpha}_{\alpha}
 & \leq & C~~.      \frac{ \vert w^{(q_n)}_{0}  \vert ^{\alpha} A_{q_{n}0}+  \vert w^{(q_{n})}_{q_{n}}  \vert ^{\alpha} A_{q_{n}q_{n}}} {\vert \vert w^{(q_{n})}\vert \vert_{\alpha}^{\alpha}}
    \end{eqnarray*}
\noindent We have that $\vert \vert w^{(q_n)} \vert \vert_{\alpha} $ goes to infinity as $n$ goes to infinity. Indeed, if not since
the sequence $\vert \vert w^{(q_n)} \vert \vert_{\alpha} $ is a increasing sequence, it must converge. Put
$w=(v_{\lambda}(i))_{i \geq 0}$ with $v_{\lambda}(0)=1$. It follows that the sequence $(w^{(q_n)})_{n \geq 0}$ converges to $w$ in
$\ell^{\alpha}$ which means that there exists a non-zero vector $w \in l^{\alpha}$ such that
$(S-\lambda I)w=0$. Hence $\lambda \in \sigma_{p}(S)$, absurd.
Now, since $ \lim A_{q_n, 0}= 0$ and $ A_{q_n, q_n}$ is bounded,  we deduce
that $\vert \vert {((S- \lambda I) u^{(q_{n})})}\vert \vert_{\alpha} $ converge to 0, and the claim is proved.
We conclude that $\lambda$ belongs to the approximate point spectrum of $S$.
$\square$

\bigskip \medskip

\noindent \textbf{Proof of Lemma \ref{lemma:left-eigenvector}:}
We introduce here another useful representation of the transition probabilities describing them column per column. Denote by
$$
\xi_m = min \{ j \ge 1 : a_j(m) \neq 0 \} \, .
$$
From \ref{tp1} and the fact $\xi_{m+1}=\zeta_m  = min \{j \geq 1 : a_{j}(m) \ne d_{j}-1\}$ , we can represent the transition probabilities in the following way: For every $m\ge 0$
\begin{equation*}
s(n,m) =
\left\{
\begin{array}{cl}
\prod_{j=1}^{\xi_m} p_j &, \ n=m-1 \, , \\
1-p_1 &, \ m=n \, , \\
(1-p_{r+1}) \prod_{j=1}^r p_j &, \ n = m + q_r -1 = m + \sum_{j=1}^{r} (d_j-1) q_{j-1}  \, , \\
& \ 1 \le r \le \xi_m -1 \, , \ \xi_m \ge 2 \, , \\
0 &, \ \textrm{otherwise} \, .
\end{array}
\right.
\end{equation*}

\smallskip

Now suppose that $v^t S = \lambda v$. Then, for all $m \geq 1$,
\begin{eqnarray}
\label{eq:rs2}
\lambda v_m & = & (vS)_m \, = \, \sum_{n=1}^\infty v_n s(n,m) \nn \\
& = & \Big( \prod_{j=1}^{\xi_m} p_j \Big) \, v_{m-1} + (1-p_1) \, v_m \nn \\
& & \qquad + \sum_{r=1}^{\xi_m -1} (1-p_{r+1}) \Big( \prod_{j=1}^r p_j \Big)
\, v_{m + q_r - 1} \, .
\end{eqnarray}
We will show later that
\begin{equation}
\label{eq:rs3}
v_{m + q_r - 1} = \frac{v_m}{ \prod_{k=1}^{r} (\iota_{\lambda}(k))^{d_k-1} } ,
\ \ 1 \le r \le \xi_m -1 \, , \ \ \xi_m \ge 2 \, .
\end{equation}
Using (\ref{eq:rs2}) and (\ref{eq:rs3}), we can use induction to prove  (\ref{eq:rs1}). Indeed, for m=1,
\begin{equation*}
\lambda v_1 = (v^t S)_1 = p_1 v_0 + (1-p_1) v_1 \quad \Rightarrow \quad
v_1 = \frac{v_0}{\iota_\lambda(1)} \, .
\end{equation*}
Now suppose that (\ref{eq:rs1}) holds for $m-1$. By (\ref{eq:rs2}) and (\ref{eq:rs3}), we have that
\begin{eqnarray}
\label{eq:rs5}
\lambda v_m & = & \Big( \prod_{j=1}^{\xi_m} p_j \Big) \, \frac{v_0}{\prod_{r=1}^{\infty} (\iota_{\lambda}(r))^{a_{r}(m-1)}} + (1-p_1) \, v_m \nn \\
& & \qquad + \sum_{r=1}^{\xi_m -1} (1-p_{r+1}) \Big( \prod_{j=1}^r p_j \Big)
\, \frac{v_m}{ \prod_{k=1}^{r} (\iota_{\lambda}(k))^{d_k-1} } \, .
\end{eqnarray}
Thus $v_0$ is equal to $v_m$ times
\begin{equation}
\label{eq:rs4}
\frac{ \prod_{r=1}^{\infty} (\iota_{\lambda}(r))^{a_{r}(m-1)} }{ \Big( \prod_{j=1}^{\xi_m} p_j \Big) } \left[ \lambda - (1-p_1) - \sum_{r=1}^{\xi_m-1} \, (1-p_{r+1}) \, \prod_{j=1}^{r}
\frac{p_j}{ \iota_\lambda(r)^{d_r-1}} \right] .
\end{equation}
Note that
$$
a_r(m-1) =
\left\{
\begin{array}{cl}
d_r - 1 &, 1 \le r < \xi_m \, , \\
a_{r}(m) - 1 &, \ r = \xi_m \, , \\
a_r(m) &, \ r \ge \zeta_m \, .
\end{array}
\right.
$$
Thus for $\zeta_m=1$ the expression in (\ref{eq:rs4}) turns out to be
$$
\prod_{r=1}^{\infty} (\iota_{\lambda}(r))^{a_{r}(m)} \, \frac{1}{\iota_\lambda(1)} \,
\left[ \frac{\lambda - (1-p_1)}{p_1} \right] = \prod_{r=1}^{\infty} (\iota_{\lambda}(r))^{a_{r}(m)} \, ,
$$
and (\ref{eq:rs1}) holds. For $\zeta_m \ge 2$, also use (\ref{eq:rs4}) and the facts
$$
\prod_{r=1}^{\infty} (\iota_{\lambda}(r))^{a_{r}(m-1)} =
\Big( \prod_{r =1}^{\xi_m-1} \iota_\lambda(r)^{d_r-1} \Big) \, \iota_{\lambda}(\xi_m))^{a_{\xi_m}(m)-1} \, \Big( \prod_{r > \xi_m} (\iota_{\lambda}(r))^{a_{r}(m)} \Big)
$$
and
$$
\begin{array}{l}
\frac{\Big( \prod_{r =k}^{\xi_m-1} \iota_\lambda(r)^{d_r-1} \Big)}{\Big( \prod_{j=k}^{\xi_m} p_j \Big)}
\left[
\iota_\lambda(k-1)^{d_{k-1}} - (1-p_k) - \sum_{r=k}^{\xi_m-1}
\, (1-p_{r+1}) \, \prod_{j=k}^{r} \frac{p_j}{ \iota_\lambda(r)^{d_r-1}}
\right] \\
=
\frac{\Big( \prod_{r =k}^{\xi_m-1} \iota_\lambda(r)^{d_r-1} \Big)}{\Big( \prod_{j=k+1}^{\xi_m} p_j \Big)}
\left[
\iota_\lambda(k) - \frac{1}{p_k} \sum_{r=k}^{\xi_m-1}
\, (1-p_{r+1}) \, \prod_{j=k}^{r} \frac{p_j}{ \iota_\lambda(r)^{d_r-1}}
\right] \\
=
\frac{\Big( \prod_{r =k+1}^{\xi_m-1} \iota_\lambda(r)^{d_r-1} \Big)}{\Big( \prod_{j=k+1}^{\xi_m} p_j \Big)}
\left[
\iota_\lambda(k)^{d_k} - (1-p_{k+1}) - \sum_{r=k+1}^{\xi_m-1}
\, (1-p_{r+1}) \, \prod_{j=k+1}^{r} \frac{p_j}{ \iota_\lambda(r)^{d_r-1}}
\right] \, .
\end{array}
$$
to verify (\ref{eq:rs1}) by finite induction.

\smallskip
We finish with the proof of (\ref{eq:rs3}). We use induction again. We should keep in mind that $\zeta_m \ge 2$. For $1 \le k \le d_1 -1$, we have
$$
\lambda v_{m+k} = p_1 v_{m+k-1} + (1-p_1)v_{m+k} \quad \Rightarrow \quad
v_{m+k} = \frac{v_m+k-1}{\iota_\lambda(1)} \, .
$$
Therefore
$$
v_{m+q_1-1} = \frac{v_m}{\iota_\lambda(1)^{d_1-1}} \, ,
$$
and (\ref{eq:rs3}) holds for every $m$ with $\xi_m \ge 2$ and $r=1$. Now fix $l\ge 1$ and suppose that (\ref{eq:rs3}) holds for every $m$ with $\xi_m \ge l+2$ and $1\le r \le l$.
For $1 \le k \le d_l -1$
\begin{eqnarray*}
\lambda v_{m+ k q_l} & = & \Big( \prod_{j=1}^{l+1} p_j \Big) \,
v_{m+\sum_{j=1}^{l-1} (d_j-1)q_{j-1} + (k-1)q_l} \nn \\
& & \qquad + (1-p_1) v_{m+kq_l} \nn \\
& & \qquad \quad + \sum_{r=1}^{l+1} (1-p_{r+1}) \Big( \prod_{j=1}^r p_j \Big)
\, v_{m+\sum_{j=1}^{r} (d_j-1)q_{j-1} + k q_l} \, ,
\end{eqnarray*}
which, by the induction hypothesis, is equal to
\begin{eqnarray*}
\lambda v_{m+ k q_l} & = & \Big( \prod_{j=1}^{l+1} p_j \Big) \,
\frac{v_{m+(k-1)q_l}}{\prod_{j=1}^{l}\iota_\lambda (j)^{d_j-1} }
+ (1-p_1) v_{m+kq_l} \nn \\
& & \qquad + \sum_{r=1}^{l+1} (1-p_{r+1}) \Big( \prod_{j=1}^r p_j \Big)
\, \frac{v_{m+ k q_l}}{\prod_{j=1}^r \iota_\lambda (j)^{d_j-1}} \, .
\end{eqnarray*}
The last expression is similar to (\ref{eq:rs5}) and wields
\begin{equation}
\label{eq:rs6}
\frac{ \prod_{r=1}^{l} (\iota_{\lambda}(r))^{d_r-1} }{ \Big( \prod_{j=1}^{l+1} p_j \Big) } \left[ \lambda - (1-p_1) - \sum_{r=1}^{l} \, (1-p_{r+1}) \, \prod_{j=1}^{r}
\frac{p_j}{ \iota_\lambda(r)^{d_r-1}} \right] .
\end{equation}
which is analogous to (\ref{eq:rs4}). Thus (\ref{eq:rs6}) is equal to $\iota_\lambda(l+1)$ and we obtain that
$$
v_{m+k q_l} = \frac{v_{m+(k-1)q_l}}{\iota_\lambda (l+1)} \, .
$$
Therefore
$$
v_{m+ (d_{l+1}-1) q_l} = \frac{v_{m}}{\iota_\lambda (l+1)^{d_{l+1} -1}} \, ,
$$
which implies that
$$
v_{m + q_{l+1} - 1} = v_{m + \sum_{k=0}^l (d_{k+1} - 1) q_k} = \frac{v_{m + (d_{l+1}-1)q_l}}{ \prod_{k=1}^{l} (\iota_{\lambda}(k))^{d_k-1} } = \frac{v_{m}}{ \prod_{k=1}^{l+1} (\iota_{\lambda}(k))^{d_k-1} } \, ,
$$
and (\ref{eq:rs3}) holds for $r=l+1$.
$\square$

\bigskip \medskip

 \noindent{\bf Acknowledgment}

\vspace{1em}

The authors would like to thank El Houcein El Abdalaoui, Sylvain Bonnot and Olivier Sester for fruitful discussions.
They also thank Christian Mauduit for fruitful discussions and also for indicating the class of adding machines based 
on Cantor Systems of numeration.


\begin{thebibliography}{99}
\bibitem{am} E. H. Abdalaoui, A. Messaoudi,  {\it On the spectrum of stochastic pertubations of the shift and Julia sets}, Fundamenta Matematicae 218(1), 47--68, 2012.
\bibitem{cg} L. Carleson, T. Gamelin, {\it Complex Dynamics}, Springer, 1993.
\bibitem{Co} J.B. Conway, {\it Functions of one Complex Variable}, Second edition, Springer Verlag, 1978.
\bibitem {DS} N. Dunford and J. T. Schwartz, {\it Linear Operators} (Interscience, 1963).
\bibitem {HW} G. H. Hardy, E. M. Wright, {\it An Introduction to the Theory of Numbers},
Oxford University Press, 1954.
\bibitem{KT1} P. Kirschenhofer, R. F. Tichy, {\it On the distribution of digits in Cantor representations of integers}, J. Number Theory 18 , 121--134, 1984.
\bibitem{KT2} P. Killeen, T.Taylor, {\it  How the Propagation of Error Through Stochastic Counters Affects Time
Discrimination and Other Psychophysical Judgments}, Psychological Review, Vol. 107(3), 430--459, 2000.
\bibitem{kt} P. Killeen, T. Taylor, {\it A stochastic adding machine and complex dynamics}, Monlinearity 13, 1998--1903, 2000.
\bibitem{l} G.F. Lawler, {\it Introduction to Stochastic Processes}, Chapman and Hall, New York, 1995.
\bibitem{ms} A. Messaoudi, D. Smania, {\it Eigenvalues of Fibonacci stochastic adding
machine}, Stochastic and Dynamics, 10(2), 291--313, 2010.
\bibitem{msv} A. Messaoudi, O. Sester, G. Valle, {\it Stochastic Adding Machines and Fibered Julia Sets}, Stochastic and Dynamics,  13(3), [26 pages], 2013.
\bibitem{Mil} J. Milnor,  Dynamics in One Complex Variable, {\it Princeton University} Press, 2006.
\bibitem{O}  J.Y. Ouvrard, {\it Probabilit\'es},  Volume 2, Cassini, Paris, 2009.
\bibitem{ses} O. Sester, {\it Hyperbolicit\'e des polyn\^omes fibr\'es}, Bull. Soc. math. France, 127, 393--428, 1999.
\bibitem{Y} K. Yosida, \it {Functional Analysis}, Springer, 1980.

\end{thebibliography}
\end{document}